\documentclass[9pt,twocolumn]{IEEEtran}
\usepackage{amsmath,amssymb,amsthm,mathrsfs,euscript,yfonts,psfrag,latexsym,dsfont,graphicx,bbm,color,amstext,wasysym,subfig,parskip,url,soul,bm,cite,balance,enumitem}
\usepackage{algorithm,algpseudocode,algorithmicx}
\usepackage{xpatch}
\makeatletter
\AtBeginDocument{\xpatchcmd{\@thm}{\thm@headpunct{.}}{\thm@headpunct{}}{}{}}
\makeatother
\usepackage[symbol]{footmisc}
\usepackage[dvipsnames]{xcolor}
\usepackage[makeroom]{cancel}
\graphicspath{{./},{./figures/}}

\DeclareMathOperator*{\argmin}{arg\:min}
\DeclareMathOperator*{\arginf}{arg\:inf}

\begin{document}
\newtheorem{thm}{Theorem}
\newtheorem{corollary}[thm]{Corollary}
\newtheorem{conj}[thm]{Conjecture}
\newtheorem{lemma}[thm]{Lemma}
\newtheorem{proposition}[thm]{Proposition}
\newtheorem{problem}{Problem}
\newtheorem{remark}{Remark}
\newtheorem{definition}{Definition}
\newtheorem{example}{Example}

\newcommand{\prox}{\rm{prox}}
\newcommand{\bp}{{\bm{p}}}
\newcommand{\bq}{{\bm{q}}}
\newcommand{\bc}{{\bm{c}}}
\newcommand{\be}{{\bm{e}}}
\newcommand{\br}{{\bm{r}}}
\newcommand{\bw}{{\bm{w}}}
\newcommand{\bx}{{\bm{x}}}
\newcommand{\by}{{\bm{y}}}
\newcommand{\bz}{{\bm{z}}}
\newcommand{\bbf}{{\bm{f}}}
\newcommand{\brh}{{\bm{\varrho}}}
\newcommand{\balpha}{{\bm{\alpha}}}
\newcommand{\bbvarphi}{{\bm{\bvarphi}}}
\newcommand{\bdelta}{{\bm{\delta}}}
\newcommand{\bPi}{{\bm{\Pi}}}
\newcommand{\bC}{{\bm{C}}}
\newcommand{\bL}{{\bm{L}}}
\newcommand{\bI}{{\bm{I}}}
\newcommand{\bP}{{\bm{P}}}
\newcommand{\bD}{{\bm{D}}}
\newcommand{\bQ}{{\bm{Q}}}
\newcommand{\epin}{{ \epsilon^{-1}}}
\renewcommand{\dagger}{{T}}
\newcommand{\ropt}{{\rm{opt}}}

\allowdisplaybreaks

\makeatletter
\newcommand{\StatexIndent}[1][3]{%
  \setlength\@tempdima{\algorithmicindent}%
  \Statex\hskip\dimexpr#1\@tempdima\relax}
\algdef{S}[WHILE]{WhileNoDo}[1]{\algorithmicwhile\ #1}%
\makeatother

\renewcommand{\thefootnote}{\fnsymbol{footnote}}

\newcommand{\differential}{{\rm{d}}}
\newcommand{\hess}{{\mathbf{Hess}}}

\newcommand{\interior}[1]{%
  {\kern0pt#1}^{\mathrm{o}}%
}

\newcommand{\bR}{\mathbb{R}}
\newcommand{\diag}{\operatorname{diag}}
\newcommand{\tr}{\operatorname{trace}}
\newcommand{\ignore}[1]{}

\newcommand{\magenta}{\color{magenta}}
\newcommand{\red}{\color{red}}
\newcommand{\blue}{\color{blue}}
\newcommand{\gray}{\color{gray}}
\definecolor{grey}{rgb}{0.6,0.3,0.3}
\definecolor{lgrey}{rgb}{0.9,.7,0.7}

\renewcommand{\qedsymbol}{\hfill\ensuremath{\blacksquare}}

\def\spacingset#1{\def\baselinestretch{#1}\small\normalsize}
\setlength{\parindent}{20pt}
\setlength{\parskip}{12pt}
\spacingset{1}

\title{\huge{Wasserstein Proximal Algorithms for the Schr{\"o}dinger Bridge Problem:\\Density Control with Nonlinear Drift}}

\author{Kenneth F. Caluya,  and Abhishek Halder
\thanks{K.F.~Caluya, and A.~Halder are with the Department of Applied Mathematics, University of California, Santa Cruz,
        CA 95064, USA;
        {\texttt{\{kcaluya,ahalder\}@ucsc.edu}}. }
\thanks{This research was partially supported by NSF award 1923278.}
}

\markboth{\today}{}

 \maketitle

\begin{abstract}
We study the Schr{\"o}dinger bridge problem (SBP) with nonlinear prior dynamics. In control-theoretic language, this is a problem of minimum effort steering of a given joint state probability density function (PDF) to another over a finite time horizon, subject to a controlled stochastic differential evolution of the state vector. As such, it can be seen as a stochastic optimal control problem in continuous time with endpoint density constraints -- a topic that originated in the physics literature in 1930s, and in the recent years, has garnered burgeoning interest in the systems-control community. 

For generic nonlinear drift, we reduce the SBP to solving a system of forward and backward Kolmogorov partial differential equations (PDEs) that are coupled through the boundary conditions, with unknowns being the ``Schr\"{o}dinger factors" -- so named since their product at any time yields the optimal controlled joint state PDF at that time. We show that if the drift is a gradient vector field, or is of mixed conservative-dissipative nature, then it is possible to transform these PDEs into a pair of initial value problems (IVPs) involving the same forward Kolmogorov operator. Combined with a recently proposed fixed point recursion that is contractive in the Hilbert metric, this opens up the possibility to numerically solve the SBPs in these cases by computing the Schr\"{o}dinger factors via a single IVP solver for the corresponding (uncontrolled) forward Kolmogorov PDE. The flows generated by such forward Kolmogorov PDEs, for the two aforementioned types of drift, in turn, enjoy gradient descent structures on the manifold of joint PDFs with respect to suitable distance functionals. We employ a proximal algorithm developed in our prior work, that exploits this geometric viewpoint, to solve these IVPs and compute the Schr\"{o}dinger factors via weighted scattered point cloud evolution in the state space. We provide the algorithmic details and illustrate the proposed framework of solving the SBPs with nonlinear prior dynamics by numerical examples. 
\end{abstract}


\section{Introduction}\label{IntroSection}
\begin{figure}[t]
\centering
\includegraphics[width=\linewidth]{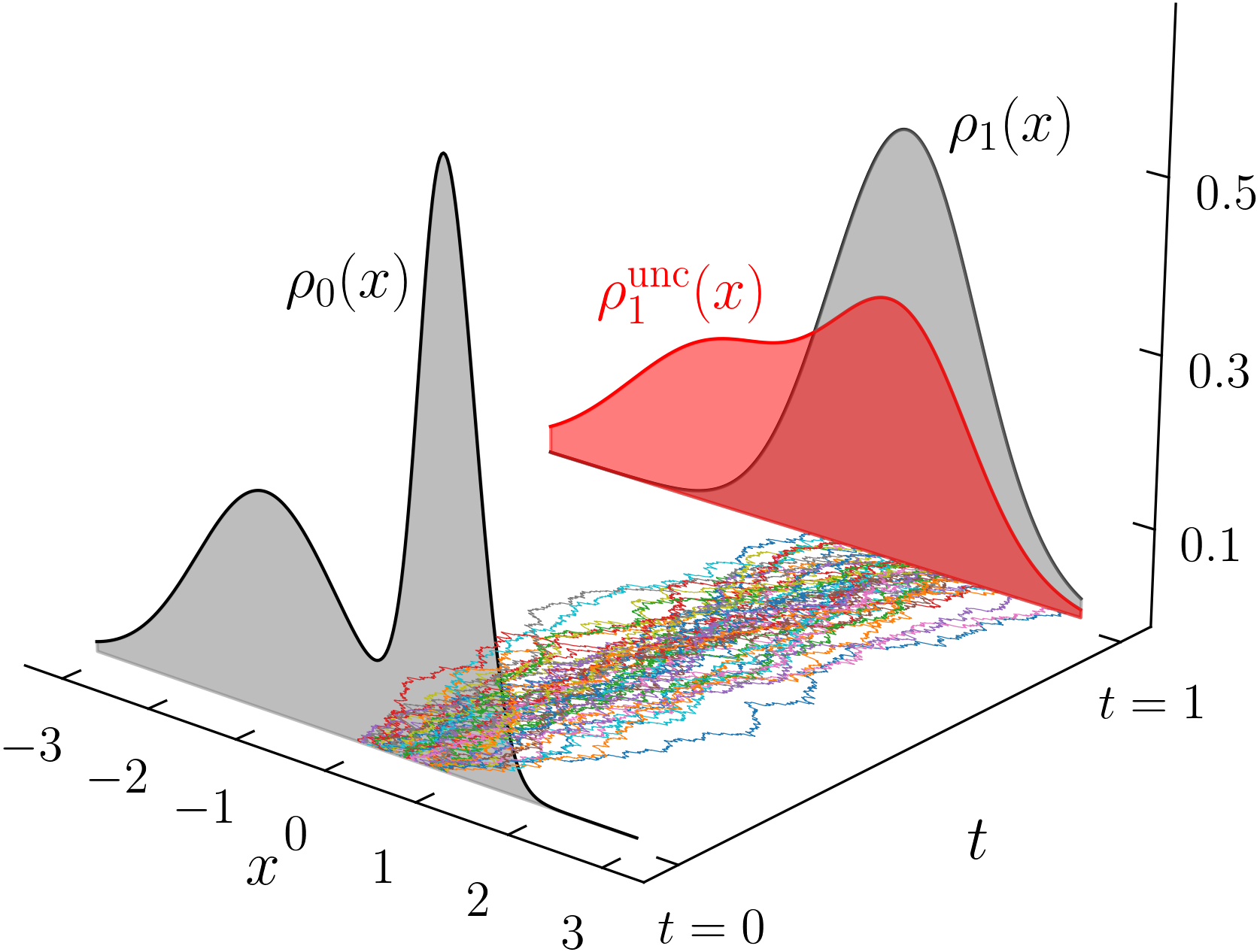}
\caption{\small{The classical SBP with 1-dimensional state space shown above, with state variable $x\in\mathbb{R}$, concerns determining the optimal control $u^{\text{opt}}(x,t)$ that steers the prescribed initial state PDF $\rho_{0}(x)$ at time $t=0$ to the prescribed terminal state PDF $\rho_{1}(x)$ at time $t=1$, while minimizing $\mathbb{E} \left\{\int_{0}^{1} \frac{1}{2} \lvert u(x,t)  \rvert^{2}  \: \differential t\right\}$ subject to a controlled diffusion (\ref{SBoriginalSDE}), i.e., the SBP solves (\ref{SBPoriginal}) wherein the expectation operator in the objective is w.r.t. the controlled state PDF $\rho(x,t)$. The optimal controlled sample paths $x^{\text{opt}}(t)$ for 100 randomly chosen initial states are shown in the $(x,t)$ plane. In the absence of control, starting from $\rho_{0}(x)$, the (uncontrolled) state PDF at $t=1$ becomes $\rho_{1}^{\text{unc}}(x)$, as depicted. For numerically solving the classical SBP with $\epsilon = 0.5$ and $\rho_{0},\rho_{1}$ as shown above, we used the fixed point recursion proposed in \cite{chen2016entropic}; see Section \ref{SubsecSBPandSchrodingerSystem} for details.}}
\label{FigSBPintro}
\vspace*{-0.1in}
\end{figure}

The Schr{\"o}dinger bridge problem (SBP) is a non-standard \emph{finite horizon} stochastic optimal control problem in continuous time. The ``non-standard" aspect stems from the fact that the SBP concerns with steering the flow of the joint state probability density function (PDF), and not the state trajectory \emph{per se}, in a controlled stochastic dynamical system 
while minimizing the expected control effort to do so. In other words, the SBP is a two-point density control problem subject to controlled trajectory-level dynamics. Fig. \ref{FigSBPintro} shows a schematic of the SBP.

Solving the SBP amounts to \emph{feedback synthesis for ensemble shaping}. The problem is of broad contemporary interest due to its potential applications in controlling a \emph{physical population} such as robotic swarm \cite{bandyopadhyay2017probabilistic}, ensemble of neurons \cite{monga2018synchronizing}, and density of highway traffic \cite{chien1997traffic}. 
 From a probabilistic perspective, one can alternatively interpret controlling the joint state PDF as that of dynamically reshaping uncertainties in a feedback loop -- a viewpoint that promotes ``control \emph{of} uncertainties" in lieu of the usual ``control \emph{with} uncertainties". That the role of feedback could be strategically synthesizing, instead of simply mitigating the uncertainties, is a recent line of thought \cite{brockett2007optimal,mazurenko2011dynamic,brockett2012notes}. 

During 1931-32, the SBP was introduced by Erwin Schr{\"o}dinger in the two articles \cite{schrodinger1931umkehrung,schrodinger1932theorie}. Schr{\"o}dinger's original motivation was to come up with a classical reformulation of quantum mechanics via diffusion processes; see e.g., \cite{zambrini1986stochastic,zambrini1986variational}. Early mathematical treatment \cite{fortet1940resolution,beurling1960automorphism,jamison1975markov} of the subject considered the case of no prior dynamics, which is what we refer to as the ``classical SBP". The stochastic optimal control interpretation for the same emerged in \cite{mikami1990variational,dai1991stochastic}. We point the readers to \cite{wakolbinger1990schrodinger,leonard2014survey} for survey of the classical SBP. More recent works \cite{chen2015stochastic,chen2015fast,chen2015optimalACC,chen2016relation,caluya2019finite,caluya2019multiinput,caluya2020reflected} have considered SBP with prior dynamics. Except the case of linear prior dynamics with Gaussian endpoint PDFs, the SBP with prior dynamics--while most relevant in control applications--remains computationally challenging in general. The computational difficulty is particularly evident in the case of \emph{nonlinear} prior dynamics, since then, first order optimality conditions for the SBP lead to solving a system of coupled nonlinear partial differential equations (PDEs). Numerical algorithms for solving the SBP are only beginning to appear now \cite{pavon2018data,li2018computations}. The purpose of this paper is to show that if the nonlinear prior dynamics is either a gradient vector field, or has both dissipative (gradient) and
conservative (Hamiltonian) components, then we can design scalable point cloud algorithms based on the recursive evaluation of certain proximal operators with respect to (w.r.t.) the Wasserstein metric, that solves the associated SBPs. The types of nonlinearities considered in this paper subsume well-known stochastic models in practical applications such as the Nyquist-Johnson resistors \cite{brockett1979stochastic,liberzon2000nonlinear}, rotational diffusions for describing the orientation in liquid crystals \cite[Ch. 16.5]{chirikjian2016harmonic}, and statistical mechanics of polymers \cite[Ch. 17]{chirikjian2016harmonic}.

The Wasserstein proximal recursions we employ, solve the forward-backward infinite-dimensional flows in the space of the \emph{Schr{\"o}dinger factors} (see Section \ref{SubsecSBPandSchrodingerSystem}), defined as functions of the state vector and time, whose product at any given time equals the optimal controlled joint state PDF at that time. These factors satisfy a system of boundary-coupled PDEs, which we solve via the proximal recursions. An intriguing aspect of our proposed framework is that the Wasserstein proximal recursions originate in the theory of optimal mass transport (OMT) \cite{villani2003topics}, and two connections between the OMT and the SBP are well-known: \emph{first}, the dynamic version \cite{benamou2000computational} of the OMT can be recovered as the zero-noise limit \cite{mikami2004monge,leonard2012schrodinger} of the classical SBP; \emph{second}, the classical SBP can be recovered as the solution of the dynamic OMT problem with Fisher information regularization \cite[Sec. 5]{chen2016relation}, \cite[Sec. 4]{conforti2017extremal}. The use of Wasserstein gradient flows to solve the SBP with prior nonlinear drift, as proposed here, offers a \emph{third} connection between the OMT and the SBP. This can be of independent interest.

The remaining of this paper is structured as follows. In Section \ref{SectionBackground}, we recap the basics of the SBP. Section \ref{SectionProblemFormulation} considers the SBP with generic nonlinear prior dynamics, deduces the necessary conditions of optimality for the same, and derives an associated Schr{\"o}dinger system via the Hopf-Cole transform. In Section \ref{SectionReformulation}, we specialize the development in Section \ref{SectionProblemFormulation} for two important types of prior nonlinear dynamics: gradient drift, and mixed conservative-dissipative drift. We show how to reformulate the corresponding Schr{\"o}dinger systems to make the same algorithmically amenable. In Section \ref{SectionProxAlgo}, we show that the reformulations derived in Section \ref{SectionReformulation} can be solved via Wasserstein proximal recursions. Section \ref{SectionNumericalExamples} provides numerical examples to illustrate the proposed algorithmic framework. Section \ref{SectionConclusions} concludes the paper.



\vspace*{-0.1in}
\section{Background}\label{SectionBackground}
\vspace*{-0.1in}
\subsubsection*{Notation}
We use bold-faced capital letters for matrices, and bold-faced lower-case letters for column vectors. The Euclidean gradient, divergence, Laplacian and Hessian operators are denoted by $\nabla$, $\nabla\cdot$, $\Delta$, and ${\mathbf{Hess}}(\cdot)$, respectively. We will sometimes put subscript to these differential operators to clarify that the operator is w.r.t. the subscripted variable (e.g, $\nabla_{\bx}$ to mean that the gradient is w.r.t. vector $\bx$); we will drop the subscript when there is no scope for confusion. We use $\langle \cdot, \cdot\rangle$ to denote the standard Euclidean inner product. For vectors $\bm{x},\bm{y} \in \mathbb{R}^n$, we have $\langle \bm{x},\bm{y} \rangle := \bm{x}^{\top}\bm{y}$, and that the squared Euclidean 2-norm $\parallel\bx\parallel_{2}^{2} \: := \langle\bx,\bx\rangle$. For matrices $\bm{A},\bm{B}$ of appropriate dimensions, $\langle \bm{A},\bm{B} \rangle:= {\mathrm{trace}}({\bm{A}}^{\top}\bm{B})$ denotes the Frobenius inner product. For matrix argument, $\|\cdot\|_{2}$ denotes the spectral norm. The symbol $\bm{0}$ denotes either a column vector or a matrix (depending on context), with all entries equal to zeros; the symbol $\bm{1}$ denotes the column vector containing all ones; $\delta(\bx-\bm{y})$ denotes the Dirac delta located at $\by$. The collection of all joint PDFs supported on $\mathbb{R}^{n}$ with finite second moments, is denoted by $\mathcal{P}_{2}\left(\mathbb{R}^{n}\right)$, i.e.,
\[\mathcal{P}_{2}\left(\mathbb{R}^{n}\right) := \bigg\{\rho : \mathbb{R}^{n} \mapsto \mathbb{R}_{\geq 0} \bigg\vert \int_{\mathbb{R}^{n}}\!\!\!\rho\:\differential\bx = 1, \int_{\mathbb{R}^{n}}\!\!\!\!\parallel\bx\parallel_{2}^{2}\rho\:\differential\bx < \infty\bigg\}.\]
We use the shorthand $\bm{x} \sim \rho$ to mean that the random vector $\bm{x}$ has the joint PDF $\rho$. The notation $\mathcal{N}\left(\bm{\mu},\bm{\Sigma}\right)$ stands for a joint Gaussian PDF with mean $\bm{\mu}$ and covariance $\bm{\Sigma}$. The symbol $\bm{I}$ denotes an identity matrix of appropriate dimension. We use the the notations $\odot$ and $\oslash$ for element-wise (i.e., Hadamard) multiplication and division, respectively. The operands $\log(\cdot)$ and $\geq$ are to be understood element-wise. The set of natural numbers is denoted by $\mathbb{N}$.

In the following, we briefly review the classical SBP, i.e., the case of no prior dynamics. We then point out a recent extension: the SBP with linear prior dynamics. Along the way, we introduce the Schr\"{o}dinger system, and the Schr\"{o}dinger factors, which will be important in our development.

\vspace*{-0.1in}
\subsection{Classical SBP and the Schr{\"o}dinger System}\label{SubsecSBPandSchrodingerSystem}
\vspace*{-0.1in}
Let the state space be $\mathbb{R}^{n}$. The \emph{classical SBP} is a stochastic optimal control problem of the form
\begin{subequations}  
\begin{align}
& \qquad\underset{\bm{u} \in \mathcal{U}} {\text{inf}}
& &  \mathbb{E} \left\{\int_{0}^{1} \frac{1}{2} \lVert\bm{u}(\bx,t)\rVert_{2}^{2}  \: \differential t \right\}, \label{SBPoriginalObj}\\
& \;  \; \; \text{subject to}
& &  \differential\bx(t) =  \bm{u}(\bx,t) \: \differential t + \sqrt{2\epsilon} \: \differential\bw(t)  \label{SBoriginalSDE},\\  
& & & \bx(t=0) \sim  \rho_{0}(\bx), \quad \bx(t=1) \sim \rho_{1}(\bx), 
\end{align}
\label{SBPoriginal}
\end{subequations}
where the set of feasible controls $\mathcal{U}$ comprises of finite energy inputs, i.e., $\mathcal{U}:=\{\bm{u} : \mathbb{R}^{n}\times [0,1] \mapsto \mathbb{R}^{n} \mid \langle\bm{u},\bm{u}\rangle < \infty\}$. Here, $\bm{w}(t)$ denotes the standard Wiener process in $\mathbb{R}^{n}$, the diffusion coefficient $\epsilon > 0$ (not necessarily small), the prescribed initial and terminal joint state PDFs are $\rho_{0}$ and $\rho_{1}$, respectively. The expectation operator in (\ref{SBPoriginalObj}) is taken w.r.t. the controlled joint state PDF $\rho(\bm{x},t)$.

The SBP (\ref{SBPoriginal}) can be transcribed into the following equivalent variational problem:
\begin{subequations}  \label{SBproblem}
\begin{align}
& \quad \underset{\rho,\bm{u}} {\inf}
& & \int_{0}^{1}\int_{\bR^{n}} \frac{1}{2} \lVert \bm{u}(\bx,t) \rVert_{2}^{2} \:\rho(\bx,t)  \: \differential\bx\:\differential t, \\
& \text{subject to}
& & \frac{\partial \rho }{\partial t} + \nabla \cdot(\rho\bm{u}) = \epsilon \Delta \rho,  \label{SBProblembc1} \\  
& & & \rho(\bm{x},0) = \rho_{0}(\bx), \quad \rho(\bm{x},1) = \rho_{1}(\bx),  \label{SBProblembc2} 
\end{align}
\end{subequations}
where (\ref{SBProblembc1}) is the Fokker-Planck or Kolmogorov's forward PDE, hereafter abbreviated as the FPK PDE, associated with the controlled stochastic differential equation (SDE) (\ref{SBoriginalSDE}). Let $\mathcal{P}(\mathbb{R}^{n})$ denote the collection of all joint PDFs supported on $\mathbb{R}^{n}$. Notice that the decision variable in problem (\ref{SBproblem}) is the pair $(\rho,\bm{u})\in\mathcal{P}(\mathbb{R}^{n})\times\mathcal{U}$.

From the first order conditions of optimality for (\ref{SBproblem}), it is easy to verify that the optimal pair $(\rho^{\ropt},\bm{u}^{\ropt})$ satisfies the following system of PDEs:  
\begin{equation}
\begin{aligned} \label{SBCondsOpt} 
&\frac{\partial \psi}{\partial t} + \frac{1}{2} \lVert \nabla \psi \rVert_{2}^{2} = -\epsilon \Delta \psi, \\ 
&\frac{\partial}{\partial t}\rho^{\ropt} + \nabla \cdot(\rho^{\ropt} \nabla \psi) = \epsilon \Delta \rho^{\ropt}, 
\end{aligned}
\end{equation}
and the optimal control $\bm{u}^{\ropt}(\bx,t) = \nabla \psi(\bx,t)$.

Furthermore, (\ref{SBCondsOpt}) can be transformed into a system of \emph{linear} PDEs via the mapping $\left(\rho^{\ropt},\psi\right)\mapsto\left(\varphi,\hat{\varphi}\right)$ given by 
\begin{subequations}
\begin{align}
\varphi(\bx,t) &= \exp \left(\frac{\psi(\bx,t)}{2\epsilon} \right), \label{HopfColevarphi}\\
\hat{\varphi}(\bx,t) &=\rho^{\ropt}(\bx,t)\exp\left(-\frac{\psi(\bx,t)}{2\epsilon}\right).\label{HopfColevarphihat} 
\end{align}
\label{HopfCole}
\end{subequations}
Specifically, the transformed variables $\left(\varphi,\hat{\varphi}\right)$ satisfy the following pair of forward-backward heat equations 
\begin{subequations}\label{ForwBackHeat}
\begin{align}
\frac{\partial \varphi}{\partial t}  &= -\epsilon \Delta \varphi, \\
\frac{\partial \hat{\varphi}}{\partial t}  &= \: \epsilon \Delta \hat{\varphi},
\end{align}
\end{subequations}
with coupled boundary conditions
\begin{subequations}\label{HeatBC}
\begin{align} 
\varphi(\bx,t=0) \hat{\varphi}(\bx,t=0) &= \rho_0(\bx),\\ 
\varphi(\bx,t=1) \hat{\varphi}(\bx,t=1) &= \rho_1(\bx).
\end{align}
\end{subequations}
To ease notation, let 
\begin{align}
\varphi_{1}(\bx):=\varphi(\bx,t=1), \quad \hat{\varphi}_{0}(\bx):=\hat{\varphi}(\bx,t=0).
\label{phi1phizerohat}
\end{align}
Then, the solution of (\ref{ForwBackHeat}) can be formally written as
\begin{subequations}  \label{HeatTransient}
\begin{align}
\varphi(\bx,t) = \int_{\mathbb{R}^n} K(t,\bx,1,\by) \varphi_1(\by) \: \differential\by,\quad t \leq 1, \\
\hat{\varphi}(\bx,t) = \int_{\mathbb{R}^n} K(0,\by,t,\bx) \hat{\varphi}_0(\by) \: \differential\by,\quad t \geq0, 
\end{align}
\end{subequations}
where 
\begin{equation} \label{HeatKernel}
K(t,\bx,s,\by) := \left(4 \pi \epsilon (t-s)\right)^{-n/2}\exp\left(-\frac{\lVert \bx-\by\rVert_{2}^2}{4\epsilon(t-s)} \right)
\end{equation}
is the heat kernel or the Markov kernel associated with the pure diffusion SDE
$\differential\bx(t) = \sqrt{2\epsilon}\:\differential\bm{w}(t)$.

Combining (\ref{HeatBC}) and (\ref{HeatTransient}), it becomes clear that finding the minimizer for the classical SBP amounts to solving for the pair $(\varphi_1,\hat{\varphi}_0)$ which satisfies the following system of nonlinear integral equations, referred to as the \emph{Schr{\"o}dinger system}, given by
\begin{subequations} \label{SchroSystem}
\begin{align}
\rho_0(\bx) = \hat{\varphi}_0(\bx) \int_{\mathbb{R}^n} K(0,\bx,1,\by) \varphi_1(\by) \: \differential\by, \\
\rho_1(\bx) = \varphi_{1}(\bx) \int_{\mathbb{R}^n} K(0,\by,1,\bx) \hat{\varphi}_0(\by) \: \differential\by .
\end{align}
\end{subequations}
The existence and uniqueness of solutions to the Schr{\"o}dinger system (\ref{SchroSystem}) were established in \cite{fortet1940resolution,jamison1975markov,beurling1960automorphism}. To compute the pair $(\varphi_1,\hat{\varphi}_0)$ from (\ref{SchroSystem}), a fixed point recursion was proposed in \cite{chen2016entropic}. Such a recursion was also proved \cite[Sec. III]{chen2016entropic} to be contractive in Hilbert's projective metric \cite{hilbert1895gerade,bushell1973hilbert}. Once the pair $(\varphi_1,\hat{\varphi}_0)$ is obtained from (\ref{SchroSystem}), then using (\ref{HeatTransient}) one can compute the pair $\left(\varphi,\hat{\varphi}\right)$. Finally, from (\ref{HopfCole}), the original decision variables $(\rho^{\ropt},\bm{u}^{\ropt})$ can be recovered via the mapping $\left(\varphi,\hat{\varphi}\right)\mapsto\left(\rho^{\ropt},\bm{u}^{\ropt}\right)$ given by
\begin{subequations}
\begin{align}
\rho^{\text{opt}}\left(\bm{x},t\right) &= \varphi\left(\bm{x},t\right)\hat{\varphi}\left(\bm{x},t\right), \label{rhooptClassicalSBP}\\
\bm{u}^{\text{opt}}\left(\bm{x},t\right) &= 2\epsilon\nabla\log\varphi(\bm{x},t). \label{uoptClassicalSBP}
\end{align}
\label{Back2OriginalVariables}		
\end{subequations}
From (\ref{rhooptClassicalSBP}), the optimal controlled joint state PDF at any time is a product of the factors $\varphi$ and $\hat{\varphi}$ at that time, and hence we refer to $(\varphi,\hat{\varphi})$ as the \emph{Schr{\"o}dinger factors}. Notice that the factors solve the boundary-coupled PDE system (\ref{ForwBackHeat})-(\ref{HeatBC}). 

The computational steps outlined in the preceding paragraph were used to solve the classical SBP in Fig. \ref{FigSBPintro} with $\epsilon = 0.5$. Using the optimal control computed from (\ref{uoptClassicalSBP}), 100 sample paths (shown in the $(x,t)$ plane in Fig. \ref{FigSBPintro}) of the optimal closed-loop SDE were simulated via the Euler-Maruyama scheme with time step-size $10^{-3}$.  

\vspace*{-0.1in}
\subsection{SBP with Linear Prior Dynamics}
\vspace*{-0.1in}
Recently, the classical SBP has been extended \cite{chen2016optimal} to the case when the prior dynamics is a linear time-varying (LTV) system, i.e., (\ref{SBoriginalSDE}) is replaced with the more general controlled SDE
\begin{equation} \label{LTVSDEcontrolled}
 \differential\bx(t) = \bm{A}(t)\bx(t) \: \differential t + \bm{B}(t)\bm{u}(\bx,t) \: \differential t + \sqrt{2\epsilon} \bm{B}(t)\: \differential\bm{w}(t),
\end{equation}
where the system matrices $\bm{A}(t)\in \mathbb{R}^{n\times n}, \bm{B}(t)\in \mathbb{R}^{n\times m}$, $m\leq n$, and the pair $(\bm{A}(t),\bm{B}(t))$ is assumed to be controllable for all $t$. When $\bm{A}(t)$ is identically zero, and $\bm{B}(t)$ is identity, then the setup reduces to the classical SBP. A \emph{Schr{\"o}dinger system} for the LTV case can be established wherein the heat kernel (\ref{HeatKernel}) is to be replaced by the Markov kernel associated with the uncontrolled SDE 
\begin{equation} \label{LTVSDE}
 \differential\bx(t) = \bm{A}(t)\bx(t) \: \differential t + \sqrt{2\epsilon} \bm{B}(t)\: \differential\bm{w}(t).
\end{equation}
We refer the readers to \cite[Sec. 4]{chen2016optimal} for the details. 

Thanks to the availability of the Markov kernel associated with (\ref{LTVSDE}), the fixed point recursion idea and the contraction results mentioned in Section \ref{SubsecSBPandSchrodingerSystem} carry through in this case, making the computation tractable as in the classical SBP. We note that in (\ref{LTVSDEcontrolled}), the noise and the control act through the same channels, i.e., the same matrix $\bm{B}(t)$ appears as the coefficient of both. This is the case we will focus in this paper. The more general case of different coefficient matrices was treated in \cite{chen2015optimal}.


\begin{figure}[t]
\centering
\includegraphics[width=0.8\linewidth]{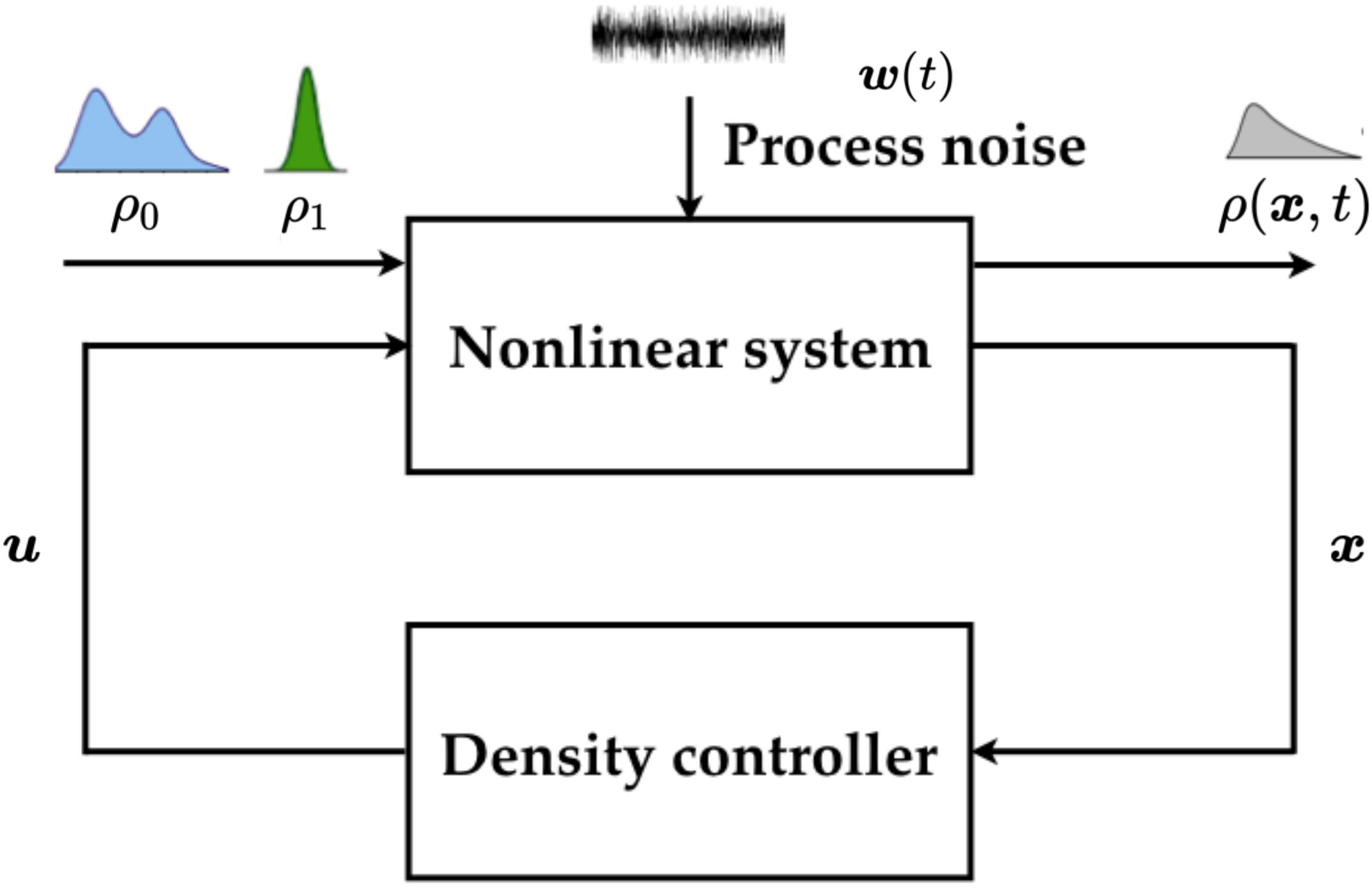}
\caption{\small{Schematic of the SBP with nonlinear prior dynamics (see Section \ref{SectionProblemFormulation} for the problem formulation).}}
\label{FigDesnityControlBlockDiagram}
\vspace*{-0.15in}
\end{figure}

\section{Problem Formulation}\label{SectionProblemFormulation}
\vspace*{-0.1in}
In this Section, we consider a much less investigated version of the SBP: a generalization of problem (\ref{SBPoriginal}) with nonlinear prior dynamics given by a (deterministic) vector field $\bm{f} : \mathbb{R}^{n} \times \mathbb{R}_{>0}\mapsto\mathbb{R}^{n}$. In particular, we address the following minimum energy stochastic optimal control problem: 
\begin{subequations}  \label{SBSDE}
\begin{align}
& \,\underset{\bm{u} \in \mathcal{U}} {\text{inf}}
& &  \mathbb{E} \left\{\int_{0}^{1} \frac{1}{2} \lVert \bm{u}(\bx,t)  \rVert_{2}^{2}  \: \differential t \right\},\\
& \text{subject to}
& &  \differential\bx = \bm{f}(\bx,t) \: \differential t +  \bm{B}(t)\bm{u}(\bx,t) \: \differential t + \sqrt{2\epsilon}  \bm{B}(t)\: \differential\bm{w}(t)  \label{SBSDEeq},\\  
& & & \bx(t=0) \sim  \rho_{0}(\bx), \quad \bx(t=1) \sim \rho_{1}(\bx),
\end{align}
\end{subequations}
where $\bx \in \mathbb{R}^{n}$, and the set $\mathcal{U}$ comprises of all finite energy inputs, as before. We make the following assumptions for $\bm{f}(\bx,t)$ and $\bm{B}(t)$:
\begin{itemize}[leftmargin=*]
\item \textbf{non-explosion and Lipschitz:} there exist constants $c_{1},c_{2}$ such that $\|\bm{f}(\bx,t)\|_{2} + \|\bm{B}(t)\|_{2} \leq c_{1}(1 + \|\bm{x}\|_{2})$, and $\|\bm{f}(\bx,t) - \bm{f}(\by,t)\|_{2} \leq c_{2}\|\bx - \by\|_{2}$ for all $\bx,\by\in\mathbb{R}^{n}$, $t\in [0,1]$;
\item \textbf{uniformly lower bounded diffusion:} there exists a constant $c_{3}$ such that the diffusion tensor $\bm{D}(t):=\bm{B}(t)\bm{B}(t)^{\top}$ satisfies $\langle\bx,\bm{D}(t)\bm{x}\rangle\geq c_{3}\|\bm{x}\|_{2}^{2}$ for all $\bm{x}\in\mathbb{R}^{n}$, $t\in [0,1]$.	
\end{itemize}
The first assumption rules out \cite[p. 66]{oksendal2013stochastic} the finite time blow up of the sample path $\bx(t)$, and ensures the existence-uniqueness for the same. The second assumption together with the first, ensures \cite[Ch. 1]{friedman2008partial} that for $1 \geq t>s\geq 0$, the transition density $K(s,\bm{y},t,\bm{x})$ associated with the uncontrolled nonlinear SDE
\begin{align}
\differential\bx(t) = \bm{f}(\bx,t)\:\differential t + \sqrt{2\epsilon}\:\bm{B}(t)\:\differential\bm{w}(t),
\label{UncontrolledNonlinearSDE}	
\end{align}
is strictly positive and everywhere continuous in $\bm{x},\bm{y}\in\mathbb{R}^{n}$. Clearly, (\ref{LTVSDEcontrolled}) corresponds to the special case $\bm{f}(\bx,t)\equiv \bm{A}(t)\bx(t)$ in (\ref{SBSDEeq}). Given $\bm{f}(\bx,t)$, $\bm{B}(t)$ and $\epsilon$, our objective is to steer the joint state PDF $\rho(\bx,t)$ from a prescribed initial PDF $\rho_{0}$ at $t=0$ to another prescribed terminal PDF $\rho_{1}$ at $t=1$ while minimizing the expected control effort. A schematic for the problem is shown in Fig. \ref{FigDesnityControlBlockDiagram}. 

Problem (\ref{SBSDE}) can be formally recast into a ``fluid dynamics'' version \cite{benamou2000computational} given by: 
\begin{subequations}  \label{SBwDrift}
\begin{align}  \label{SBwDriftCostFunctional}
& \qquad \underset{(\rho,\bm{u})} {\text{inf}} 
& &  \frac{1}{2} \int_{0}^{1}\int_{\bR^n}  \lVert \bm{u}(\bx,t) \rVert_{2}^{2} \:\rho(\bx,t)  \: \differential\bx \: \differential t  \\ 
& \text{subject to}
& & \frac{\partial \rho }{\partial t} + \nabla \cdot\left(\rho \left(\bm{f}+\bm{B}(t)\bm{u}\right)\right) = \epsilon \langle\bm{D}(t),\hess\left(\rho\right)\rangle,\label{SBFluidsFP} \\  
& & & \rho(\bx,0) = \rho_{0}(\bx), \quad \rho(\bx,1) = \rho_{1}(\bx), \label{SBFluidsBC}
\end{align}
\end{subequations}
where the infimum is taken over all pairs $(\rho,\bm{u})\in\mathcal{P}(\mathbb{R}^{n})\times\mathcal{U}$ satisfying  (\ref{SBFluidsFP})-(\ref{SBFluidsBC}). We note that (\ref{SBFluidsFP}) is the controlled FPK PDE which governs the flow of the joint PDF associated with the SDE (\ref{SBSDEeq}). Conceptually, problem (\ref{SBwDrift}) is to problem (\ref{SBSDE}) what problem (\ref{SBproblem}) is to problem (\ref{SBPoriginal}).

%

\vspace*{-0.1in}
\subsection{Existence and Uniqueness} 
\vspace*{-0.1in}
The change of variable $(\rho,\bm{u}) \mapsto (\rho,\bm{m})$ given by $\bm{m} := \rho\bm{u}$, results in the following convex reformulation of (\ref{SBwDrift}):
\begin{subequations}\label{SBmrho}
\begin{align} 
& \underset{(\rho,\bm{m})} {\text{minimize}} 
\quad \frac{1}{2}\int_{0}^{1}\int_{\bR^n}  J(\rho,\bm{m})\: \differential\bx\:\differential t  \\
& \text{subject to}
\;\:\label{SBmrhobc1} \frac{\partial \rho }{\partial t} + \nabla \cdot\left(\rho\bm{f} + \bm{B}(t)\bm{m}\right) - \epsilon\langle\bm{D}(t),\hess\left(\rho\right)\rangle = 0 ,  \\  
& \qquad\qquad\;\label{SBmrhobc2} \rho(\bx,0) = \rho_{0}(\bx), \quad \rho(\bx,1) = \rho_{1}(\bx),  
\end{align}
\end{subequations}
where \[J(\rho,\bm{m}):= \begin{cases}\parallel\bm{m}\parallel_{2}^{2}/\rho & \text{if}\quad \rho>0,\\
 0 & \text{if}\quad(\bm{m},\rho)=(\bm{0},0),\\
 +\infty &\text{otherwise}.	
 \end{cases}\] 
 Since $J(\rho,\bm{m})$ is the perspective function of the strictly convex map $\bm{m}\mapsto \lVert\bm{m}\rVert_{2}^{2}$, hence $J$ is jointly strictly convex in $(\rho,\bm{m})$. The constraints (\ref{SBmrhobc1})-(\ref{SBmrhobc2}) are linear in $(\rho,\bm{m})$. However, proving the existence of minimizing pair $(\rho^{{\rm{opt}}},\bm{m}^{{\rm{opt}}})$, or equivalently $(\rho^{{\rm{opt}}},\bm{u}^{{\rm{opt}}})$, is nontrivial. It turns out that the existence-uniqueness for the solution of (\ref{SBwDrift}) or (\ref{SBmrho}), can be guaranteed \cite[Theorem 3.2]{jamison1974reciprocal},\cite[Sec. 10]{wakolbinger1990schrodinger} for compactly supported $\rho_0,\rho_1$, provided the transition density for (\ref{UncontrolledNonlinearSDE}) is strictly positive and everywhere continuous in $\bm{x},\bm{y}\in\mathbb{R}^{n}$. This is where the aforesaid assumptions on $\bm{f}(\bm{x},t),\bm{B}(t)$ come into play. The details for two specific cases of interest, viz. gradient and mixed conservative-dissipative drift, are given in Appendix \ref{AppendinxRegularity}.

\vspace*{-0.1in} 
\subsection{Conditions for Optimality}
\vspace*{-0.1in}
 In this Section, we show that the first order optimality conditions for (\ref{SBwDrift}) correspond to a coupled system of nonlinear PDEs. Consider the Lagrangian associated with (\ref{SBwDrift}):
\begin{equation}
\begin{aligned}\label{SBwDriftLagrangian}
&\mathcal{L}(\rho,\bm{u},\psi) :=  \int_{0}^{1}\!\!\int_{\bR^n} \bigg\{\frac{1}{2}\lVert \bm{u}(\bx,t) \rVert_{2}^{2}\:\rho(\bx,t) + \psi(\bx,t) \times \bigg(  \frac{\partial \rho}{\partial t}\\
& + \nabla \cdot \left(\left(\bm{f}+\bm{B}(t)\bm{u}\right)\rho(\bx,t)\right) - \epsilon\langle\bm{D}(t),\hess\left(\rho\right)\rangle\bigg) \bigg \}
\: \differential\bx\:\differential t,
\end{aligned}
\end{equation}
where $\psi(\bx,t)$ is a $C^{1}\left(\mathbb{R}^{n};\mathbb{R}_{>0}\right)$ Lagrange multiplier. Let 
\begin{equation}
\begin{aligned}
\mathcal{P}_{01}\left(\mathbb{R}^{n}\right) := \bigg \{ \rho(\bx,t) \mid\: &\rho \geq 0,  \int_{\mathbb{R}^{n}} \!\!\rho\:\differential\bx =1, \\ 
&\rho(\bx,0)= \rho_0,\rho(\bx,1)=\rho_1 \bigg  \}.
 \end{aligned}
\end{equation}Performing the unconstrained minimization of the Lagrangian $\mathcal{L}$  over $\mathcal{P}_{01}\left(\mathbb{R}^{n}\right)\times\mathcal{U}$ yields the following result (proof in Appendix \ref{AppendixCondOpt}).

\begin{proposition}\label{propSBPnonlineargeneral} (\textbf{Optimal control and optimal state PDF})
The pair $(\rho^{{\rm{opt}}}(\bx,t),\bm{u}^{{\rm{opt}}}(\bx,t))$ that solves (\ref{SBwDrift}), must satisfy the system of coupled PDEs
\begin{subequations} \label{SBwDriftOptConds}
\begin{align} 
&\frac{\partial \psi }{\partial t } + \frac{1}{2}\lVert \bm{B}(t)^{\!\top}\nabla\psi \rVert_{2}^{2} + \langle\nabla\psi,\bm{f}\rangle = - \epsilon \langle\bm{D}(t),\hess\left(\psi\right)\rangle, \label{HJBgenNonlinear}\\ 
     &\frac{\partial}{\partial t}\rho^{{\rm{opt}}} + \nabla\!\cdot\!(\rho^{{\rm{opt}}}(\bm{f} + \bm{B}(t)^{\!\top}\nabla \psi)) = \epsilon\langle\bm{D}(t),\hess\left(\rho^{{\rm{opt}}}\right)\rangle, \label{FPKgenNonlinear}
     \end{align} 
\end{subequations}\\
with boundary conditions  
\begin{equation}
\rho^{{\rm{opt}}}(\bx,0)=\rho_0(\bx),\quad \rho^{{\rm{opt}}}(\bx,1) = \rho_1(\bx),\label{SBwDriftBC}
\end{equation}
 and \[\bm{u}^{{\rm{opt}}}(\bx,t) = \bm{B}(t)^{\top}\nabla\psi(\bx,t).\]
\end{proposition}

\begin{remark} Notice that the system of coupled PDEs (\ref{SBwDriftOptConds}) is the same as that appearing in classical mean field games, but instead of the usual boundary conditions (see e.g., \cite[equation (2)]{lasry2007mean}) 
\begin{equation}
\rho^{{\rm{opt}}}(\bx,0) = \rho_0(\bx)\;\text{(given)}, \quad \psi(\bx,1) = \psi_1(\bx)\;\text{(given)},
\end{equation}
our boundary conditions are given by (\ref{SBwDriftBC}).
\end{remark}
The PDE (\ref{HJBgenNonlinear}) is the Hamilton-Jacobi-Bellman (HJB) equation while the PDE (\ref{FPKgenNonlinear}) is the controlled FPK equation. Computing the optimal pair $(\rho^{\ropt},\bm{u}^{\ropt})$, or equivalently $(\rho^{\ropt},\psi)$ form (\ref{SBwDriftOptConds})-(\ref{SBwDriftBC}) calls for solving a system of coupled nonlinear PDEs with atypical boundary conditions, and is challenging in general. 

To tackle the issues of coupling and nonlinearity, we resort to the \emph{Hopf-Cole transform} \cite{hopf1950partial,cole1951quasi} given by (\ref{HopfCole}). 
The geometric interpretation of this transform for optimal control problems were investigated in \cite{leger2019hopf,leger2018geometric}. Similar ideas have also appeared in the stochastic control literature \cite{fleming1982logarithmic,dai1990markov,bakshi2020schrodinger}. In our context, this transform allows to convert the HJB-FPK system of coupled nonlinear PDEs (\ref{SBwDriftOptConds})-(\ref{SBwDriftBC}) into a system of \emph{boundary-coupled linear PDEs}. We summarize this in the following Theorem (proof in Appendix \ref{AppendixHopfColeProof}).

\begin{thm}\label{ThmHopfCole}
(\textbf{Hopf-Cole transform})
Given $\bm{f}$, $\epsilon$, $\rho_{0}$, $\rho_{1}$, consider the Hopf-Cole transform $(\rho^{{\rm{opt}}},\psi)\mapsto(\varphi,\hat{\varphi})$ defined via (\ref{HopfCole}), applied to the system of nonlinear PDEs (\ref{SBwDriftOptConds})-(\ref{SBwDriftBC}). Then the pair $(\varphi,\hat{\varphi})$ satisfies the system of linear PDEs
\begin{subequations}  \label{FPKBK}
    \begin{align}
          \frac{\partial \varphi }{\partial t} & = - \langle\nabla\varphi,\bm{f}\rangle -\epsilon\langle\bm{D}(t),\hess\left(\varphi\right)\rangle, \label{BK} \\ 
         \frac{\partial \hat{\varphi}}{\partial t} &=  - \nabla \cdot( \hat{\varphi}\bm{f}) + \epsilon \langle\bm{D}(t),\hess\left(\hat{\varphi}\right)\rangle,\label{FPK}
         \end{align}
\end{subequations}
with boundary conditions (borrowing notations (\ref{phi1phizerohat}))
\begin{equation} \label{FPKBKbcs}
\varphi_{0}(\bx)\hat{\varphi}_{0}(\bx) = \rho_{0}(\bx),\quad \varphi_{1}(\bx) \hat{\varphi}_{1}(\bx) = \rho_{1}(\bx).
\end{equation} 
Moreover, the optimal controlled state PDF $\rho^{{\rm{opt}}}(\bx,t)$  solving (\ref{SBwDrift}) is given by (\ref{rhooptClassicalSBP}). The optimal control for the same is given by
\[\bm{u}^{{\rm{opt}}}(\bx,t) = 2\epsilon\bm{B}(t)^{\top}\nabla\log\varphi.\]
\end{thm}

We note that (\ref{BK}) is the \emph{backward Kolmogorov equation} in variable $\varphi$, and (\ref{FPK}) is the \emph{forward Kolmogorov} or the FPK equation in variable $\hat{\varphi}$, associated with (\ref{UncontrolledNonlinearSDE}). Indeed, (\ref{FPKBK}) generalizes the forward-backward heat equations (\ref{ForwBackHeat}). As expected, setting $\bm{f}\equiv\bm{0}$ and $\bm{B}(t)\equiv\bm{I}$ in (\ref{FPKBK}) recovers (\ref{ForwBackHeat}). Following the nomenclature in Section \ref{SubsecSBPandSchrodingerSystem}, the pair $(\varphi,\hat{\varphi})$ solving (\ref{FPKBK})-(\ref{FPKBKbcs}) defines the \emph{Schr\"{o}dinger factors} for problem (\ref{SBwDrift}).

The essence of Theorem \ref{ThmHopfCole} is that instead of solving the system of coupled nonlinear PDEs (\ref{SBwDriftOptConds}), we can solve the system of linear PDEs (\ref{FPKBK}) provided that we  compute the pair $(\varphi_1,\hat{\varphi}_0)$ which serves as the endpoint data for 
\begin{subequations}  \label{FPKB2}
    \begin{align} 
          \frac{\partial \varphi }{\partial t} & = - \langle\nabla\varphi,\bm{f}\rangle -\epsilon\langle\bm{D}(t),\hess\left(\varphi\right)\rangle,
         \, \varphi(\bx,1) = \varphi_1(\bx),\\ 
         \frac{\partial \hat{\varphi}}{\partial t} &=  - \nabla \cdot( \hat{\varphi}\bm{f}) + \epsilon \langle\bm{D}(t),\hess\left(\hat{\varphi}\right)\rangle,\, \hat{\varphi}(\bx,0) = \hat{\varphi}_0(\bx). 
         \end{align}
\end{subequations}
Denoting the forward and backward Kolmogorov operators in (\ref{FPKB2}) as $L_{\rm{FK}}$ and $L_{\rm{BK}}$ respectively, we can write (\ref{FPKB2}) succinctly as an infinite dimensional two point boundary value problem
\begin{align}\label{FPKB2vec} 
          \frac{\partial}{\partial t}\begin{pmatrix}\varphi\\
	\hat{\varphi}
\end{pmatrix}
 = \begin{pmatrix}
 L_{\rm{BK}} & 0\\
 0 & L_{\rm{FK}}	
 \end{pmatrix}\begin{pmatrix}\varphi\\
	\hat{\varphi}
\end{pmatrix}, \: \begin{pmatrix}\varphi(\bx,1)\\
 	\hat{\varphi}(\bx,0)
 \end{pmatrix}
 = \begin{pmatrix}\varphi_{1}(\bx)\\  
 \hat{\varphi}_0(\bx)
\end{pmatrix}. 
\end{align}
Solving (\ref{FPKB2}) or (\ref{FPKB2vec}) will yield $(\varphi(\bx,t),\hat{\varphi}(\bx,t))$ for all $t\in [0,1]$, which in turn can be used to determine $ (\rho^{\ropt}(\bx,t),\bm{u}^{\ropt}(\bx,t))$ via (\ref{Back2OriginalVariables}). Notice that this computational pipeline hinges on the ability to first compute the pair $(\varphi_1,\hat{\varphi}_0)$ using (\ref{FPKBKbcs}) \emph{and} (\ref{FPKB2}), and then utilize $(\varphi_1,\hat{\varphi}_0)$ to solve (\ref{FPKB2}) or (\ref{FPKB2vec}). However, unlike the situation in Section \ref{SectionBackground}, the difficulty now is that the closed-form expression of the Markov kernel associated with (\ref{UncontrolledNonlinearSDE}) is not available in general. This prevents us from setting up a \emph{Schr\"{o}dinger system} like (\ref{SchroSystem}) to solve for the pair $(\varphi_1,\hat{\varphi}_0)$.

In the following Section \ref{SectionReformulation}, we will reformulate the Schr\"{o}dinger system for two cases: when $\bm{f}$ is gradient of a potential, and when the prior drift has mixed conservative-dissipative structure, i.e., a \emph{degenerate} diffusion of the form (\ref{SBSDEeq}) with $\bm{u,w}\in\mathbb{R}^{m}$, the state space dimension $n \equiv 2m$, and a constant input matrix $\bm{B} \in \mathbb{R}^{n\times m}$. In Section \ref{SectionProxAlgo}, we will then show how these reformulations can harness the recently-introduced proximal recursions \cite{caluya2018proximal,caluya2019gradient} for computing the pair $(\varphi_1,\hat{\varphi}_0)$, and consequently the pair $(\varphi(\bx,t),\hat{\varphi}(\bx,t))$. This is appealing since these proximal algorithms do not require spatial discretization or function approximation, and instead evolve weighted scattered point cloud data. Therefore, such algorithms hold the promise for solving the SBP with nonlinear prior dynamics in high dimensions when no analytical handle on the Markov kernel is available.

\vspace*{-0.1in}
\section{Reformulation of the Schr\"{o}dinger Systems}\label{SectionReformulation}
\vspace*{-0.1in}
Next, we provide algorithmically tractable reformulations of the Schr\"{o}dinger systems for two important types of prior nonlinear dynamics: gradient drift in Section \ref{SubsectionSchrodingerSystemwithGradientDrift}, and mixed conservative-dissipative drift (i.e., degenerate diffusion, see e.g., \cite[Sec. 7, 8]{brockett1997notes}) in Section \ref{SubsectionSchrodingerSystemwithMixedDrift}. For clarity of exposition, in Section \ref{SubsectionSchrodingerSystemwithGradientDrift}, we consider $\bm{B}(t)\equiv\bm{I}$. In Section \ref{SubsectionSchrodingerSystemwithMixedDrift}, the matrix $\bm{B}$ will be non-identity.

\subsection{The Case of Gradient Drift}\label{SubsectionSchrodingerSystemwithGradientDrift}
\vspace*{-0.1in}
We now consider $\bm{f}$ to be a gradient vector field, i.e., $\bm{f}= -\nabla V(\bx)$ for some $C^{2}\left(\mathbb{R}^{n}\right)$ function $V:\mathbb{R}^{n} \mapsto \mathbb{R}_{\geq 0}$, and $\bm{B}(t)=\bm{I}$. This reduces (\ref{SBwDrift}) to the following:
\begin{subequations}  \label{SBwGradDrift}
\begin{align} 
& \qquad\quad \underset{(\rho,\bm{u})} {\text{inf}}
& &  \frac{1}{2} \int_{0}^{1}\int_{\bR^n}  \lVert \bm{u}(\bx,t) \rVert_{2}^{2} \:\rho(\bx,t)  \: \differential\bx \: \differential t , \label{SBwGradDriftObj}\\ 
&\; \; \; \text{subject to}
& & \frac{\partial \rho }{\partial t} + \nabla \cdot(\rho (\bm{u}-\nabla V)) = \epsilon \Delta \rho, \\  
& & & \rho(\bx,0) = \rho_{0}(\bx), \quad \rho(\bx,1) = \rho_{1}(\bx).
\end{align}
\end{subequations}
Applying Proposition \ref{propSBPnonlineargeneral} and Theorem \ref{ThmHopfCole} to (\ref{SBwGradDrift}), we arrive at the system of linear PDEs for $(\varphi,\hat{\varphi})$ given by
\begin{subequations}\label{SchrodingerSystemGradient}
\begin{align}
          \frac{\partial \varphi }{\partial t} & =   \langle\nabla\varphi,\nabla V\rangle  -\epsilon \Delta \varphi , \label{GradBK}\\ 
         \frac{\partial \hat{\varphi}} {\partial t} &=  \nabla \cdot( \hat{\varphi}\nabla V) + \epsilon \Delta  \hat{\varphi},\label{GradFK}
 \end{align}
\end{subequations}
with boundary conditions 
\begin{equation}\label{BCforSBPgrad}
\varphi_0(\bx)\hat{\varphi}_0(\bx) = \rho_0(\bx), \quad  \varphi_1(\bx) \hat{\varphi}_1(\bx) = \rho_1(\bx).
\end{equation}
As expected, (\ref{GradBK}) is a backward Kolmogorov PDE, and (\ref{GradFK}) is a forward Kolmogorov or FPK\footnote{This particular instance of the FPK PDE (\ref{GradFK}) is also known as the ``Smoluchowski equation" \cite[Ch. II.4.(vi)]{chandrasekhar1943stochastic}.} PDE.
 
We would like to exploit the structure of (\ref{SchrodingerSystemGradient}) to develop an algorithm that computes the pair $(\varphi_1,\hat{\varphi}_0)$ which can then serve as the terminal and initial data for the following system solving for the pair $(\varphi(\bx,t),\hat{\varphi}(\bx,t))$:
\begin{subequations} \label{KolmGradDrift}
    \begin{align}
          \frac{\partial \varphi }{\partial t} & =  \langle\nabla\varphi,\nabla V\rangle  -\epsilon \Delta \varphi \label{BKgradrift},
         \quad \varphi(\bx,1) = \varphi_1(\bx), \\ 
         \frac{\partial \hat{\varphi}}{\partial t} &=   \nabla \cdot( \hat{\varphi}\nabla V) + \epsilon \Delta  \hat{\varphi},\quad \hat{\varphi}(\bx,0) = \hat{\varphi}_0(\bx). \label{FPKgradrift}
         \end{align}
\end{subequations}
After having $(\varphi(\bx,t),\hat{\varphi}(\bx,t))$ from (\ref{KolmGradDrift}), we can obtain the pair $(\rho^{\ropt},\bm{u}^{\ropt})$ via (\ref{Back2OriginalVariables}). To this end, the following result (proof in Appendix \ref{AppendixProofGradTVP2IVPTheorem}) will be an important step.

\begin{thm}\label{GradTVP2IVPTheorem}
Given $V(\bx), \epsilon, \varphi_{1}(\bx)$, consider the terminal value problem (TVP) (\ref{BKgradrift}) in unknown $\varphi(\bx,t)$.
Let $s:=1-t$, and $q(\bx,s) := \varphi(\bx,t) = \varphi(\bx,1-s)$. Then $q$ satisfies the initial value problem (IVP) 
\begin{equation}
\frac{\partial q}{\partial s} = -\langle\nabla q,\nabla V\rangle + \epsilon \Delta  q, \quad q(\bx,0) = \varphi_1(\bx)  \label{IVP1}.
\end{equation}
Further, $p(\bx,s):= q(\bx,s)\exp\left(-V(\bx)/\epsilon\right)$ satisfies the IVP 
\begin{equation} \label{IVP2}
\frac{\partial p}{\partial s} = \nabla \cdot (p \nabla V) + \epsilon \Delta p, \quad p(\bx,0) = q(\bx,0)\exp\left(-V(\bx)/\epsilon\right),
\end{equation}
where $q$ is a smooth solution of (\ref{IVP1}), and $V$ is such that 
\[\int_{\mathbb{R}^n}q(\bx,0)\exp\left(-V(\bx)/\epsilon \right) \: \differential\bx < \infty, \quad \text{for all}\;\epsilon > 0.\]
\end{thm}
\noindent Thanks to Theorem \ref{GradTVP2IVPTheorem}, solving (\ref{KolmGradDrift}) amounts to solving the system
\begin{subequations} \label{transPDE}
\begin{align} 
        \frac{\partial \hat{\varphi}}{\partial t}& =   \nabla \cdot( \hat{\varphi}\nabla V) + \epsilon \Delta  \hat{\varphi} ,\quad \hat{\varphi}(\bx,0) = \hat{\varphi}_0(\bx),   \label{transPDE1} \\
        \frac{\partial p}{\partial s} &= \nabla \cdot (p \nabla V) + \epsilon \Delta p , \quad p(\bx,0) = \varphi_{1}(\bx)\exp\left(-V(\bx)/\epsilon \right). \label{transPDE2}
 \end{align}
\end{subequations}
Notice that (\ref{transPDE1}) and (\ref{transPDE2}) involve the exact same PDE with different initial conditions, to be integrated in different time coordinates $t$ and $s$, where $t=1-s$. This implies that availability of a single FPK IVP solver is enough to set up a fixed point recursion for the pair $(\varphi_{1},\hat{\varphi}_{0})$ via (\ref{transPDE}). Once $p$ is computed from (\ref{transPDE2}), we can recover $\varphi$ by the relation $\varphi(\bx,t)=\varphi(\bx,1-s)=p(\bx,s)/\exp\left(-V(\bx)/\epsilon \right)$.

\subsection{The Case of Mixed Conservative-Dissipative Drift}\label{SubsectionSchrodingerSystemwithMixedDrift}
\vspace*{-0.1in}
We now consider a \emph{degenerate} diffusion of the form (\ref{SBSDEeq}) with $\bm{u,w}\in\mathbb{R}^{m}$. The state $\bx$ consists of two sub-vectors $\bm{\xi},\bm{\eta}\in\mathbb{R}^{m}$, i.e., $\bx := \left(\bm{\xi}, \bm{\eta}\right)^{\top}\in\mathbb{R}^{n}$ with $n \equiv 2m$. The controlled SDE is given by
\begin{align}
\begin{pmatrix}
\differential\bm{\xi}\\
\differential\bm{\eta}	
\end{pmatrix} &= \bigg\{\underbrace{\begin{pmatrix}
 \bm{\eta}\\
 -\nabla_{\bm{\xi}}V\left(\bm{\xi}\right) - \kappa\bm{\eta}	
 \end{pmatrix}}_{\bm{f}(\bx)} + \underbrace{\begin{pmatrix} 
 \bm{0}_{m\times m}\\
 \bm{I}_{m\times m}	
 \end{pmatrix}}_{\bm{B}}
\bm{u}\bigg\}\differential t \nonumber\\
&\quad+ \sqrt{2\epsilon\kappa}\underbrace{\begin{pmatrix} 
 \bm{0}_{m\times m}\\
 \bm{I}_{m\times m}	
 \end{pmatrix}}_{\bm{B}}
\differential\bm{w}, \quad \kappa > 0,
\label{DegenSDE}	
\end{align}
i.e., we consider the SBP (\ref{SBSDE}) with the constraint (\ref{SBSDEeq}) replaced by (\ref{DegenSDE}). Here, we assume that $V\in C^{2}(\mathbb{R}^{m})$, $\inf V > -\infty$, and that $\hess\left(V\right)$ is uniformly lower bounded.
To glean physical motivation, one may consider $m=3$ and think of $\bm{\xi}$ and $\bm{\eta}$ as three-dimensional position and velocity vectors, respectively \cite[Sec. V.B]{caluya2019gradient}; see also \cite[Examples 2,3]{liberzon2000nonlinear}. So, (\ref{DegenSDE}) is a controlled Langevin equation. We are thus led to solve an instance of (\ref{SBwDrift}) wherein for the constraint (\ref{SBFluidsFP}), we use the vector field $\bm{f}$ and the matrix $\bm{B}$ as shown above, and set the diffusion tensor as $\bm{D}:=\kappa\bm{B}\bm{B}^{\top}\in\mathbb{R}^{n\times n}$.

Using Proposition \ref{propSBPnonlineargeneral}, the first order conditions of optimality for this variant of SBP yields the coupled HJB-FPK system:   
\begin{subequations}\label{DegenrateCoupledFPKHJB}
\begin{align}
\frac{\partial \psi }{\partial t} & =   -\langle \bm{\eta},\nabla _{\bm{\xi}}\psi  \rangle + \langle \nabla_{\bm{\xi}}V\left(\bm{\xi}\right) + \kappa\bm{\eta},  \nabla_{\bm{\eta}}\psi\rangle  -\epsilon \kappa \Delta_{\bm{\eta}} \psi \nonumber\\
&\qquad - \frac{1}{2} \parallel\nabla_{\bm{\eta}} \psi\parallel_{2}^{2},\\
\frac{\partial \rho^{{\rm{opt}}} }{\partial t} &= -\langle \bm{\eta},\nabla _{\bm{\xi}}\rho^{{\rm{opt}}}  \rangle + \nabla_{\bm{\eta}} \cdot\left(\rho^{{\rm{opt}}}\left(\nabla_{\bm{\xi}}V\left(\bm{\xi}\right) + \kappa\bm{\eta} - \nabla_{\bm{\eta}}\psi\right)\right) \nonumber\\
&\qquad+ \epsilon \kappa \Delta_{\bm{\eta}} \rho^{{\rm{opt}}},	
\end{align}
\end{subequations}
with boundary conditions $\rho^{{\rm{opt}}}\left(\bx,0\right) = \rho_{0}(\bx)$ and $\rho^{{\rm{opt}}}\left(\bx,1\right) = \rho_{1}(\bx)$, and the optimal control $\bm{u}^{{\rm{opt}}}(\bx,t) = \nabla_{\bm{\eta}}\psi\left(\bx,t\right)$.

Following Theorem \ref{ThmHopfCole}, we now apply the Hopf-Cole transform $\left(\psi,\rho^{{\rm{opt}}}\right) \mapsto \left(\varphi,\hat{\varphi}\right)$ given by (\ref{HopfCole}), to the system (\ref{DegenrateCoupledFPKHJB}). This results in the following system of linear PDEs for $(\varphi,\hat{\varphi})$:
\begin{subequations}\label{SchrodingerSystemDegen}
\begin{align}
          \frac{\partial \varphi }{\partial t} & =   -\langle \bm{\eta},\nabla _{\bm{\xi}}\varphi  \rangle + \langle \nabla_{\bm{\xi}}V\left(\bm{\xi}\right) + \kappa\bm{\eta},  \nabla_{\bm{\eta}}\varphi\rangle  -\epsilon \kappa \Delta_{\bm{\eta}} \varphi ,\label{BKdegen}\\ 
         \frac{\partial \hat{\varphi}} {\partial t} &= -\langle \bm{\eta},\nabla _{\bm{\xi}}\hat{\varphi}  \rangle + \nabla_{\bm{\eta}} \cdot\left(\hat{\varphi}\left(\nabla_{\bm{\xi}}V\left(\bm{\xi}\right) + \kappa\bm{\eta}\right)\right) + \epsilon \kappa \Delta_{\bm{\eta}} \hat{\varphi}, \label{FKdegen}
 \end{align}
\end{subequations}
with boundary conditions (\ref{BCforSBPgrad}). In particular, (\ref{BKdegen}) is a backward Kolmogorov PDE, and (\ref{FKdegen}) is a forward Kolmogorov or FPK\footnote{This particular instance of the FPK PDE (\ref{FKdegen}) is also known as the ``kinetic Fokker-Planck equation" \cite{villani2006hypocoercive}.} PDE.

Next, we establish a result (proof in Appendix \ref{AppendinxDegenTVP2IVPTheorem}) for (\ref{SchrodingerSystemDegen}) that is similar in flavor to Theorem \ref{GradTVP2IVPTheorem}, and will be useful for designing proximal algorithm in Section \ref{SectionProxAlgo}.

\begin{thm}\label{DegenTVP2IVPTheorem}
Given $\epsilon,\varphi_{1}({\bm{\xi},\bm{\eta}})$, consider the TVP comprising of the PDE (\ref{BKdegen}) in variable $\varphi$ together with the boundary condition $\varphi(\bm{\xi},\bm{\eta},1)=\varphi_{1}({\bm{\xi},\bm{\eta}})$. Let $s:=1-t$, and 
$q(\bm{\xi},\bm{\eta},s):= \varphi(\bm{\xi},\bm{\eta},t) =\varphi(\bm{\xi},\bm{\eta},1-s)$. Then $q$ satisfies the IVP
\begin{subequations}
\begin{align}
    &\frac{\partial q }{\partial s} =   \langle \bm{\eta},\nabla _{\bm{\xi}}q  \rangle -\langle \nabla_{\bm{\xi}}V\left(\bm{\xi}\right) + \kappa\bm{\eta},  \nabla_{\bm{\eta}}q\rangle  +\epsilon \kappa \: \Delta_{\bm{\eta}} q, \\ 
    &q(\bm{\xi},\bm{\eta},0) =\varphi_{1}({\bm{\xi},\bm{\eta}}). 
\end{align}    
\label{BKdegenIVP}	
\end{subequations}
Further, let $\bm{\vartheta}:=-\bm{\eta}$, and 
\[\widetilde{p}(\bm{\xi},-\bm{\eta},s) := q(\bm{\xi},\bm{\eta},s) \exp\!\left(\!-\frac{1}{\epsilon} \left(\!\frac{1}{2} \lVert \bm{\eta}\rVert_{2}^{2} + V(\bm{\xi})\!\right)\!\right).\]
Then, $p(\bm{\xi},\bm{\vartheta},s) := \widetilde{p}(\bm{\xi},-\bm{\eta},s)$ satisfies the IVP
\begin{subequations}
\begin{align}
 &\frac{\partial p} {\partial s} = -\langle \bm{\vartheta},\nabla _{\bm{\xi}}p  \rangle + \nabla_{\bm{\vartheta}} \cdot\left(p \left(\nabla_{\bm{\xi}}V\left(\bm{\xi}\right) + \kappa\bm{\vartheta}\right)\right) + \epsilon \kappa\:\Delta_{\bm{\vartheta}} p,\\
&p(\bm{\xi},\bm{\vartheta},0)= q(\bm{\xi},\bm{\eta},0) \exp\left( -\frac{1}{\epsilon} \left( \frac{1}{2} \lVert \bm{\eta}\rVert_{2}^{2} + V(\bm{\xi}) \right) \right),
\end{align}
\label{FKdegenIVP}	
\end{subequations}
where $q$ is a smooth solution of (\ref{BKdegenIVP}), and $V$ is such that 
\[\int_{\mathbb{R}^{m}} \int_{\mathbb{R}^{m	}}q(\bm{\xi},\bm{\eta},0) \exp\left( -\frac{1}{\epsilon} \left( \frac{1}{2} \lVert \bm{\eta}\rVert_{2}^{2} + V(\bm{\xi}) \right) \right) \: \differential \bm{\xi}\: \differential \bm{\eta} < \infty,\]
for all $\epsilon >0$.
\end{thm}

\noindent Thanks to Theorem \ref{DegenTVP2IVPTheorem}, solving the system (\ref{SchrodingerSystemDegen}) with boundary conditions
\[\varphi(\bm{\xi},\bm{\eta},1) = \varphi_{1}(\bm{\xi},\bm{\eta}), \quad \hat{\varphi}(\bm{\xi},\bm{\eta},0) = \hat{\varphi}_{0}(\bm{\xi},\bm{\eta}),\]
reduces to solving the system
\begin{subequations} \label{transPDEdegen}
\begin{align} 
        \frac{\partial \hat{\varphi}} {\partial t} &= -\langle \bm{\eta},\nabla _{\bm{\xi}}\hat{\varphi}  \rangle + \nabla_{\bm{\eta}} \cdot\left(\hat{\varphi}\left(\nabla_{\bm{\xi}}V\left(\bm{\xi}\right) + \kappa\bm{\eta}\right)\right) + \epsilon \kappa \Delta_{\bm{\eta}} \hat{\varphi},   \label{transPDE1degen} \\
        \frac{\partial p} {\partial s} &= -\langle \bm{\vartheta},\nabla _{\bm{\xi}}p  \rangle + \nabla_{\bm{\vartheta}} \cdot\left(p \left(\nabla_{\bm{\xi}}V\left(\bm{\xi}\right) + \kappa\bm{\vartheta}\right)\right) + \epsilon \kappa\:\Delta_{\bm{\vartheta}} p, \label{transPDE2degen}
 \end{align}
\end{subequations}
with initial conditions
\begin{subequations} \label{ivpPDEdegen}
\begin{align} 
\hat{\varphi}(\bm{\xi},\bm{\eta},0) &= \hat{\varphi}_{0}(\bm{\xi},\bm{\eta}), \label{ivpPDE1degen}\\
p(\bm{\xi},\bm{\vartheta},0) &= q(\bm{\xi},\bm{\eta},0) \exp\left( -\frac{1}{\epsilon} \left( \frac{1}{2} \lVert \bm{\eta}\rVert_{2}^{2} + V(\bm{\xi}) \right) \right)\nonumber\\
&= \varphi_{1}(\bm{\xi},-\bm{\vartheta}) \exp\left( -\frac{1}{\epsilon} \left( \frac{1}{2} \lVert \bm{\vartheta}\rVert_{2}^{2} + V(\bm{\xi}) \right) \right). \label{ivpPDE2degen}	
\end{align}
\end{subequations}
Similar to Section \ref{SubsectionSchrodingerSystemwithGradientDrift}, notice that (\ref{transPDE1degen}) and (\ref{transPDE2degen}) involve the exact same kinetic Fokker-Planck PDE with different initial conditions, to be integrated in different time coordinates $t$ and $s$, where $t=1-s$. This implies that availability of a single kinetic FPK IVP solver is enough to set up a fixed point recursion for the pair $(\varphi_{1},\hat{\varphi}_{0})$ via (\ref{transPDEdegen}). Once $p(\bm{\xi},\bm{\vartheta},s)$ is computed from (\ref{transPDE2degen}), we can recover $\varphi(\bm{\xi},\bm{\eta},t)$ using the relation $\varphi(\bm{\xi},\bm{\eta},t)=\varphi(\bm{\xi},\bm{\eta},1-s)=p(\bm{\xi},\bm{\vartheta},s)/\exp\left( -\frac{1}{\epsilon} \left( \frac{1}{2} \lVert \bm{\eta}\rVert_{2}^{2} + V(\bm{\xi}) \right) \right)$.

\begin{remark}\label{ReamrkDifferenceBetnGradANDDegenIVPtransfrom}
While both Theorem \ref{GradTVP2IVPTheorem} and Theorem \ref{DegenTVP2IVPTheorem} reduce the respective Schr\"{o}dinger systems--which are two point forward-backward Kolmogorov systems--to forward-forward Kolmogorov IVPs, it is important to recognize the subtle difference between the transformations employed in these two theorems.	Unlike (\ref{SchrodingerSystemGradient}), the Schr\"{o}dinger system (\ref{SchrodingerSystemDegen}) is not reversible (see e.g., \cite[Ch. 6.1]{pavliotis2014stochastic}). To account this, Theorem \ref{DegenTVP2IVPTheorem} used a spatial transform $(\bm{\xi},\bm{\eta})\mapsto(\bm{\xi},\bm{\vartheta})$, in addition to the time transform $t\mapsto s$. This resulted in a sequence of change-of-variables $\varphi \mapsto q \mapsto \widetilde{p} \mapsto p$ in Theorem \ref{DegenTVP2IVPTheorem}. To put the matter in perspective, Theorem \ref{GradTVP2IVPTheorem} involved the change-of-variables $\varphi \mapsto q \mapsto p$. 
\end{remark}

\section{Wasserstein Proximal Algorithms}\label{SectionProxAlgo}

In this Section, we will design variational recursions to solve the SBP with gradient or mixed conservative-dissipative drift. For the case of the gradient drift, the variational recursions will iteratively solve the system of IVPs (\ref{transPDE}). For the case of the mixed conservative-dissipative drift, the variational recursions will iteratively solve the system of IVPs (\ref{transPDEdegen})-(\ref{ivpPDEdegen}). The main idea here is to design certain Lyapunov functionals (usually referred to as the `free energy functionals') such that the FPK PDEs in (\ref{transPDE}) and (\ref{transPDEdegen}) can be recast as the gradient flow of the respective Lyapunov functionals w.r.t. suitable metric on the manifold $\mathcal{P}_{2}\left(\mathbb{R}^{n}\right)$, which is the space of all PDFs supported on $\mathbb{R}^{n}$ with finite second moment. For theoretical developments along this line, we refer the readers to \cite{ambrosio2008gradient}. Our recent works \cite{caluya2018proximal,caluya2019gradient} showed that such infinite dimensional gradient flow structure can be exploited for designing grid-less recursive algorithms to solve for the transient solutions of the FPK PDE. We collect the key conceptual ingredients next.

\subsection{Proximal Operators for Infinite Dimensional Gradient Flows}
First, we introduce the 2-Wasserstein metric $W$, that will play an important role in the development that follows.

\begin{definition} \label{DefWasserstein}(\textbf{2-Wasserstein metric}) The 2-Wasserstein metric between two probability measures $\mu_0,\mu_1$, each supported on $\mathbb{R}^n$, is 
\begin{equation}
 W(\mu_0,\mu_1):= \!\left(\underset{m \in \mathcal{M} \left(\mu_0,\mu_1 \right)} \inf\!  \int_{\mathbb{R}^n \times \mathbb{R}^n} \!  \lVert \bx-\bm{y} \rVert_{2}^{2} \: \differential m(\bx,\by)\right)^{\!\!\frac{1}{2}},
 \label{EqnWassDef}
\end{equation}
where $\mathcal{M} \left(\mu_0,\mu_1 \right)$ is the set of all joint probability measures (i.e., couplings) supported on the product space $\mathbb{R}^n \times \mathbb{R}^n$ whose marginals are $\mu_0$ and $\mu_1$, respectively. Whenever $\mu_0$ and $\mu_1$ are absolutely continuous, their respective PDFs $\rho_0$ and $\rho_1$ exist (i.e., $\differential\mu_{0} = \rho_{0}(\bx)\differential\bx$, $\differential\mu_{1} = \rho_{1}(\bm{y})\differential\bm{y}$), and we use the equivalent notation $W(\rho_0,\rho_1)$. 
\end{definition}
\noindent Although Definition \ref{DefWasserstein} introduces the 2-Wasserstein metric $W$ as a static variational problem, Brenier and Benamou \cite{benamou2000computational} pointed out that  $W^{2}$ admits the following equivalent dynamic reformulation: 
\begin{subequations}\label{BBformula}
\begin{align} 
\!\!\!\!\!\!\!W^2(\rho_0,\rho_1) = 
\underset{\rho,\bm{u}} {\inf}\qquad&\int_{0}^{1}\!\!\int_{\bR^{n}} \frac{1}{2} \lVert \bm{u}(\bx,t) \rVert_{2}^{2} \: \rho(\bx,t)  \: \differential\bx \differential t, \\
\text{subject to} & \: \frac{\partial \rho }{\partial t} + \nabla \cdot(\rho \bm{u}) = 0,\\  
& \: \rho(\bx,0) = \rho_{0}(\bx), \quad \rho(\bx,1) = \rho_{1}(\bx), 
\end{align}
\end{subequations}
which can be interpreted as the minimum energy control problem to steer the PDF $\rho_0$ to $\rho_1$ over the time interval $[0,1]$ subject to the deterministic dynamics $\dot{\bx}=\bm{u}$. Formally, the classical SBP (\ref{SBproblem}) can be viewed as a dynamic stochastic regularization of (\ref{BBformula}), and this regularization is captured by the diffusion term $\epsilon \Delta \rho $ in (\ref{SBProblembc1}). These ideas can be made rigorous by proving weak convergence of measures \cite{leonard2013survey,leonard2012schrodinger,mikami2004monge,mikami2008optimal}  in the $\epsilon \downarrow 0$ limit, i.e., the solution of (\ref{SBproblem}) converges to that of (\ref{BBformula}) in the said limit.

Now the idea is to compute the flows generated by the IVPs (\ref{transPDE}) or (\ref{transPDEdegen}), via variational recursions over discrete time index pair $(t_{k-1},s_{k-1}):=((k-1)\tau,(k-1)\sigma)$, $k\in\mathbb{N}$,  given by the map
\begin{align}
\begin{pmatrix}
\hat{\phi}_{t_{k-1}}\\
\varpi_{s_{k-1}}
\end{pmatrix} \mapsto \begin{pmatrix}
\hat{\phi}_{t_{k}}\\
\varpi_{s_{k}}
\end{pmatrix} := \begin{pmatrix}
\underset{\hat{\phi} \in \mathcal{P}_2({\mathbb{R}^n})} \arginf \:  \frac{1}{2}d^{2}\left(\hat{\phi}_{t_{k-1}},\hat{\phi}\right) + \tau F(\hat{\phi})\\
\underset{\varpi \in \mathcal{P}_2({\mathbb{R}^n})} \arginf \:  \frac{1}{2}d^{2}\left(\varpi_{s_{k-1}},\varpi\right) + \sigma F(\varpi)
\end{pmatrix},
\label{IVPrecursion}
\end{align}
where $\tau,\sigma>0$ are time-step sizes, and the functionals $\left(d,F\right)$ are chosen such that 
\begin{subequations}
\begin{align}
\hat{\phi}_{t_{k-1}}(\bx) &\rightarrow \hat{\varphi}(\bx,t=(k-1)\tau)&\text{in $L^{1}(\mathbb{R}^n)$ as}\quad\tau \downarrow 0,\label{Convergencephihat}\\
\varpi_{s_{k-1}}(\bx) &\rightarrow p(\bx,s=(k-1)\sigma)&\text{in $L^{1}(\mathbb{R}^n)$ as}\quad\sigma \downarrow 0.\label{Convergencep}
\end{align}
\label{ConvergenceGuarantee}
\end{subequations}
In words, the variational recursions (\ref{IVPrecursion}) are designed to approximate the flows generated by (\ref{transPDE}) or (\ref{transPDEdegen}) in the small time-step limit.

The variational maps appearing in (\ref{IVPrecursion}) resemble \emph{proximal operators} \cite{parikh2014proximal,bauschke2011convex} well-known in the Euclidean and general Hilbert space optimization literature. In our setting too, $d$ will be a distance metric on $\mathcal{P}_{2}\left(\mathbb{R}^{n}\right)$, and $F$ will be a Lyapunov functional, i.e., non-negative and decreasing along the flows generated by the IVPs (\ref{transPDE}) or (\ref{transPDEdegen}), and thus motivate defining the proximal operators  
\begin{subequations}
\begin{align}
{\rm{prox}}_{\tau F}^{d}\left(\hat{\phi}_{t_{k-1}}\right) &:= \underset{\hat{\phi} \in \mathcal{P}_2({\mathbb{R}^n})} \arginf \:  \frac{1}{2}d^{2}\left(\hat{\phi}_{t_{k-1}},\hat{\phi}\right) + \tau F(\hat{\phi}), \label{proxphihat}\\
{\rm{prox}}_{\sigma F}^{d}\left(\varpi_{s_{k-1}}\right) &:= \underset{\varpi \in \mathcal{P}_2({\mathbb{R}^n})} \arginf \:  \frac{1}{2}d^{2}\left(\varpi_{s_{k-1}},\varpi\right) + \sigma F(\varpi),
\end{align}
\label{DefProx}
\end{subequations}
which respectively read as the proximal operator of $\tau F$ and $\sigma F$, w.r.t. the metric $d$. This metric viewpoint allows clear geometric interpretation of the proximal recursions (\ref{IVPrecursion}): they define gradient descent of the functionals $\tau F$ and $\sigma F$, measured w.r.t. the distance metric $d$; see e.g., \cite{halder2017gradient,halder2018gradient,caluya2019gradient,halder2019proximal}. For (\ref{transPDE}), the distance $d$ will turn out to be the 2-Wasserstein metric $W$. For (\ref{transPDEdegen}), the distance metric $d$ will be a variant of $W$.

\subsection{Proximal Recursions} 
We now outline how the proximal recursions (\ref{DefProx}) can be constructed for the two nonlinear drifts of our interest: gradient and mixed conservative-dissipative.

\subsubsection{Gradient drift}\label{subsubsecGrad}
The seminal paper \cite{jordan1998variational} showed that the flows generated by the FPK PDEs of the form (\ref{transPDE}) can be seen as the gradient descent of the Lyapunov functional 
\begin{equation}
F_{\text{gradient}}(\cdot): =  \int_{\mathbb{R}^n} V(\bx) (\cdot) \: \differential\bx + 
\epsilon \int_{\mathbb{R}^n} \left(\cdot\right) \log\left(\cdot\right) \: \differential\bx, 
\label{GradFreeEnergy}
\end{equation}
w.r.t. the distance metric $W$ in $\mathcal{P}_2({\mathbb{R}^n})$. Here, $(\cdot)$ is a placeholder for either $\hat{\varphi}$ or $p$ in (\ref{transPDE}). From (\ref{IVPrecursion}) and (\ref{DefProx}), the solutions of the IVPs (\ref{transPDE}), denoted by $\left(\hat{\varphi}(\bx,t),p(\bx,s)\right)$, can then be approximated (as in (\ref{ConvergenceGuarantee})) by the following proximal recursions:
\begin{align}
\begin{pmatrix}
\hat{\phi}_{t_{k}}\\
\varpi_{s_{k}}
\end{pmatrix} = \begin{pmatrix}
{\rm{prox}}_{\tau F_{\text{gradient}}}^{W}\left(\hat{\phi}_{t_{k-1}}\right)\\
{\rm{prox}}_{\sigma F_{\text{gradient}}}^{W}\left(\varpi_{s_{k-1}}\right)
\end{pmatrix}, \quad k\in\mathbb{N},
\label{GradProx}
\end{align}
where $F_{\text{gradient}}$ is given by (\ref{GradFreeEnergy}), and the initial conditions: 
\begin{align*}
\hat{\phi}_{t_{0}} &=  \hat{\varphi}(\bx,0) \stackrel{\text{(\ref{transPDE1})}}{=} \hat{\varphi}_{0}(\bx),\\
\varpi_{s_{0}} &= p(\bx,0) \stackrel{\text{(\ref{transPDE2})}}{=} \varphi_{1}(\bx)\exp(-V(\bx)/\epsilon).	
\end{align*}
In other words, for the IVPs (\ref{transPDE}), the pair $(d,F) \equiv (W,F_{\text{gradient}})$.

To numerically implement the proximal recursions (\ref{GradProx}), we will use the weighted scattered point cloud algorithm developed in our previous work \cite[Sec. III]{caluya2019gradient}. Since the metric $W$ (see Definition \ref{DefWasserstein}) is itself defined as a variational problem, evaluating proximal operators w.r.t. $W$, requires solving nested functional minimization problems. The algorithm in \cite[Sec. III.B]{caluya2019gradient} uses an entropic or Sinkhorn regularization followed by block coordinate ascent in the dual space, and provides convergence guarantees via contraction mapping \cite[Sec. III.C]{caluya2019gradient}.


\subsubsection{Mixed conservative-dissipative drift}\label{subsubsecMixed}
For FPK PDEs of the form (\ref{transPDEdegen}), recall that $\bm{x}:=\left(\bm{\xi},\bm{\eta}\right)^{\top}\in\mathbb{R}^{n}$, $\bm{y}:=\left(\bm{\xi},\bm{\vartheta}\right)^{\top}\in\mathbb{R}^{n}$, and let the ``Hamiltonian-like" function $H\left(\bx\right) := \frac{1}{2}\|\bm{\eta}\|_{2}^{2} + V(\bm{\xi}) = \frac{1}{2}\|\bm{\vartheta}\|_{2}^{2} + V(\bm{\xi}) =: H(\bm{y})$. It can be verified that
\[F_{\text{mixed}}^{\text{a}}\left(\hat{\varphi}\right):= \int_{\mathbb{R}^{n}}\!\! H(\bx)\hat{\varphi}\left(\bx,t\right)\:\differential\bx + \epsilon\int_{\mathbb{R}^{n}}\!\!\hat{\varphi}\left(\bx,t\right)\log\hat{\varphi}\left(\bx,t\right)\:\differential\bx,\]
and 
\[F_{\text{mixed}}^{\text{b}}\left(p\right):= \int_{\mathbb{R}^{n}}\!\! H(\by)p\left(\by,s\right)\:\differential\by + \epsilon\int_{\mathbb{R}^{n}}\!\!p\left(\by,s\right)\log p\left(\by,s\right)\:\differential\by,\]
are Lyapunov functionals along the flows of (\ref{transPDE1degen}) and (\ref{transPDE2degen}), respectively. However, since the degenerate diffusions do not allow statistical reversibility, the flows for (\ref{transPDEdegen}) cannot be recast as the gradient descent of the above functionals w.r.t. $W$. To circumvent this, consider instead the functionals
\begin{subequations}
\begin{align}
\widetilde{F}_{\text{mixed}}^{\text{a}}\left(\hat{\varphi}\right) := \!\int_{\mathbb{R}^{n}} \!\frac{1}{2}\|\bm{\eta}\|_{2}^{2}\hat{\varphi}\left(\bx,t\right)\:\differential\bx + \epsilon\!\!\int_{\mathbb{R}^{n}}\!\!\hat{\varphi}\left(\bx,t\right)\log\hat{\varphi}\left(\bx,t\right)\:\differential\bx,\label{MixedFreeEnergya}\\
\widetilde{F}_{\text{mixed}}^{\text{b}}\left(p\right) := \!\int_{\mathbb{R}^{n}} \!\frac{1}{2}\|\bm{\vartheta}\|_{2}^{2}p\left(\by,s\right)\:\differential\by + \epsilon\!\!\int_{\mathbb{R}^{n}}\!\!p\left(\by,s\right)\log p\left(\by,s\right)\:\differential\by.\label{MixedFreeEnergyb}
\end{align}
\label{MixedFreeEnergy}		
\end{subequations}
Also, consider the following distance functionals which are modified versions of (\ref{EqnWassDef}):
\begin{subequations} 
\begin{align}
\widetilde{W}_{\tau}(\mu_0,\mu_1):= \!\left(\underset{m \in \mathcal{M} \left(\mu_0,\mu_1 \right)} \inf\!  \int_{\mathbb{R}^n \times \mathbb{R}^n} \!\!\!\! s_{\tau}\left(\bx,\overline{\bx}\right) \: \differential m(\bx,\overline{\bx})\!\!\right)^{\!\!\frac{1}{2}},\label{Wtildetau}\\
\widetilde{W}_{\sigma}(\mu_0,\mu_1):= \!\left(\underset{m \in \mathcal{M} \left(\mu_0,\mu_1 \right)} \inf\!  \int_{\mathbb{R}^n \times \mathbb{R}^n} \!\!\!\! s_{\sigma}\left(\by,\overline{\by}\right) \: \differential m(\by,\overline{\by})\!\!\right)^{\!\!\frac{1}{2}},\label{Wtildesigma}
\end{align}
 \label{Wtilde}
\end{subequations} 
where $\bx = (\bm{\xi},\bm{\eta})^{\top}$ and $\overline{\bx} = (\overline{\bm{\xi}},\overline{\vphantom{\bm{\xi}}\bm{\eta}})^{\top}$ are two realizations of the state vector in (\ref{DegenSDE}); similarly, $\by = (\bm{\xi},\bm{\vartheta})^{\top}$, $\overline{\by} = (\overline{\bm{\xi}},\overline{\vphantom{\bm{\xi}}\bm{\vartheta}})^{\top}$, and
\begin{subequations} 
\begin{align}
\!\!\!\!s_{\tau}\left(\bx,\overline{\bx}\right) := \| \overline{\vphantom{\bm{\xi}}\bm{\eta}} - \bm{\eta} + \tau\nabla V(\bm{\xi}) \|_{2}^{2} + 12 \bigg\| \frac{\overline{\bm{\xi}} - \bm{\xi}}{\tau} - \frac{\overline{\vphantom{\bm{\xi}}\bm{\eta}} + \bm{\eta}}{2} \bigg\|_{2}^{2},\label{Defstau}\\
\!\!\!\!s_{\sigma}\left(\by,\overline{\by}\right) := \| \overline{\vphantom{\bm{\xi}}\bm{\vartheta}} - \bm{\vartheta} + \sigma\nabla V(\bm{\xi}) \|_{2}^{2} + 12 \bigg\| \frac{\overline{\bm{\xi}} - \bm{\xi}}{\sigma} - \frac{\overline{\vphantom{\bm{\xi}}\bm{\vartheta}} + \bm{\vartheta}}{2} \bigg\|_{2}^{2}.\label{Defssigma}
\end{align}
\label{Defs}	
\end{subequations}
Following \cite[Scheme 2b]{duong2014conservative}, we then set up the proximal recursions
\begin{align}
\begin{pmatrix}
\hat{\phi}_{t_{k}}\\
\varpi_{s_{k}}
\end{pmatrix} = \begin{pmatrix}
{\rm{prox}}_{\kappa\tau \widetilde{F}_{\text{mixed}}^{\text{a}}}^{\widetilde{W}_{\tau}}\left(\hat{\phi}_{t_{k-1}}\right)\\
{\rm{prox}}_{\kappa\sigma \widetilde{F}_{\text{mixed}}^{\text{b}}}^{\widetilde{W}_{\sigma}}\left(\varpi_{s_{k-1}}\right)
\end{pmatrix}, \quad k\in\mathbb{N},
\label{MixedProx}
\end{align}
with initial conditions (\ref{ivpPDEdegen}). For the above recursions, the results from \cite[Theorem 2.4]{duong2014conservative} provide the consistency guarantees (\ref{ConvergenceGuarantee}) for the flows generated by (\ref{transPDEdegen}).

The weighted scattered point cloud algorithm from \cite{caluya2019gradient} that we mentioned in Section \ref{subsubsecGrad}, can be applied to this case too (see \cite[Section V.B]{caluya2019gradient}), and we will use the same to solve (\ref{transPDEdegen}) by numerically performing the proximal recursion (\ref{MixedProx}).

\subsection{Sinkhorn Proximal Algorithms for Solving (\ref{GradProx}) and (\ref{MixedProx})}\label{subsecWassProxSchrodinger}

Our development so far have reduced solving the SBPs with gradient and mixed conservative-dissipative drifts to that of solving the proximal recursions (\ref{GradProx}) and (\ref{MixedProx}), respectively. For numerical computation, we perform these proximal recursions over weighted scattered point clouds with sample size $N$. Specifically, let $\bm{x}^{i}(t)$ be the state vector for the $i$-th sample at time $t$, and for $i=1,\hdots,N$, $k\in\mathbb{N}$, define the $N\times 1$ vectors 
\begin{subequations}
\begin{align}
\bm{\hat{\phi}}_{k-1}^{i} &:= \hat{\phi}_{t_{k-1}}\left(\bm{x}^{i}(t_{k-1})\right), \label{vechatphi}\\
\bm{\varpi}_{k-1}^{i} &:= \varpi_{s_{k-1}}\left(\bm{x}^{i}(s_{k-1})\right), \label{vecvarpi}
\end{align}
\label{discreteFunc}	
\end{subequations}
wherein the superscript $i$ in the left-hand-side denotes the ``$i$-th component of". The recursions (\ref{GradProx}) and (\ref{MixedProx}) are performed over the point clouds $\left\{ \bx^{i}(t_{k-1}),\bm{\hat{\phi}}_{k-1}^{i}  \right\}_{i=1}^N$ and $\left\{ \bx^{i}(s_{k-1}),\bm{\varpi}_{k-1}^{i}  \right\}_{i=1}^N$, respectively, following the algorithm given in \cite{caluya2018proximal,caluya2019gradient}. To do so, the state vectors $\bm{x}^{i}$ are updated by applying the Euler-Maruyama\footnote{Here, we use the Euler-Maruyama scheme for ease of implementation; it can be replaced by any SDE integrator; see \cite[Remark 1 in Sec. III.B]{caluya2019gradient}.} scheme to the appropriate version of the \emph{uncontrolled} SDE (\ref{UncontrolledNonlinearSDE}). For instance, in the gradient drift case, the uncontrolled SDE associated with (\ref{transPDE}) is
 \begin{align}
\differential\bm{x} = -\nabla V(\bm{x}) \: \differential t + \sqrt{2\epsilon} \: \differential\bm{w}.
\label{GradUncontrolledSDE}
\end{align}

For $k\in\mathbb{N}$, we suppose that $\{\bx^{i}(t_{k-1})\}_{i=1}^{N}\mapsto \{\bx^{i}(t_{k})\}_{i=1}^{N}$ and $\{\bx^{i}(s_{k-1})\}_{i=1}^{N}\mapsto \{\bx^{i}(s_{k})\}_{i=1}^{N}$, are the Euler-Maruyama updates associated with (\ref{GradUncontrolledSDE}) with time step-sizes $\tau$ and $\sigma$, respectively. Following \cite{caluya2019gradient}, we write the Sinkhorn regularized proximal recursions for (\ref{GradProx}) in vector form:
{\small{\begin{subequations}
	\begin{align}
&\bm{\hat{\phi}}_k =\: \underset{\bm{\hat{\phi}}} \argmin  \bigg\{ \underset{\bm{M} \in  \Pi(\bm{\hat{\phi}}_{k-1}, \bm{\hat{\phi}})}  \min \frac{1}{2} \bigg\langle\!\bm{C}\left(\{\bx^{i}(t_{k-1}),\bx^{i}(t_{k})\}_{i=1}^{N}\right),\bm{M}\!\bigg\rangle \nonumber\\
&+ \gamma \langle \bm{M}, \log \bm{M} \rangle +  \tau \langle  V\left(\{\bm{x}^{i}(t_{k-1})\}_{i=1}^{N}\right)+ \epsilon\log\bm{\hat{\phi}},\bm{\hat{\phi}} \rangle \bigg \},\label{GradSinkhornProxphihat}\\
&\bm{\varpi}_k =\: \underset{\bm{\varpi}} \argmin   \bigg\{ \underset{\bm{M} \in  \Pi(\bm{\varpi}_{k-1}, \bm{\varpi})}  \min \frac{1}{2} \bigg\langle\! \bm{C}\left(\{\bx^{i}(s_{k-1}),\bx^{i}(s_{k})\}_{i=1}^{N}\right),\bm{M}\!\bigg\rangle \nonumber\\
&+ \gamma \langle \bm{M}, \log \bm{M} \rangle + \sigma \langle  V\left(\{\bm{x}^{i}(s_{k-1})\}_{i=1}^{N}\right)+ \epsilon\log\bm{\varpi},\bm{\varpi} \rangle \bigg \},	\label{GradSinkhornProxvarpi}
\end{align}
\label{GradSinkhornProx}
\end{subequations}}}
wherein for $k\in\mathbb{N}$, the $(i,j)$\textsuperscript{th} component of the matrix $\bm{C}\in\mathbb{R}^{N\times N}$ in (\ref{GradSinkhornProxphihat}) equals $\|\bx^{i}(t_{k-1}) - \bx^{j}(t_{k})\|_{2}^{2}$; likewise, the $(i,j)$\textsuperscript{th} component of the matrix $\bm{C}$ in (\ref{GradSinkhornProxvarpi}) equals $\|\bx^{i}(s_{k-1}) - \bx^{j}(s_{k})\|_{2}^{2}$. In (\ref{GradSinkhornProx}), the notation $\Pi(\bm{a},\bm{b})$ stands for the set of all matrices $\bm{M}\in\mathbb{R}^{N\times N}$ such that $\bm{M}\geq 0$, $\bm{M}\bm{1}=\bm{a}$, and $\bm{M}^{\top}\bm{1}=\bm{b}$, for given admissible\footnote{Here, ``admissible" means that the vectors $\bm{a},\bm{b}$ are component-wise nonnegative, and have equal element-wise sum. This is due to the fact that the generators in (\ref{transPDE}) and (\ref{transPDEdegen}) are both integral preserving and nonnegativity preserving.} vectors $\bm{a},\bm{b}\in\mathbb{R}^{N}$. Furthermore, $V\left(\{\bm{x}^{i}(t_{k-1})\}_{i=1}^{N}\right)$ returns the $N\times 1$ vector whose $i$\textsuperscript{th} element equals $V\left(\bm{x}^{i}(t_{k-1})\right)$; likewise for $V\left(\{\bm{x}^{i}(s_{k-1})\}_{i=1}^{N}\right)$. The term $\langle \bm{M}, \log \bm{M} \rangle$ is the Sinkhorn/entropic regularization, and  $\gamma>0$ is a small regularization parameter. Dualizing (\ref{GradSinkhornProx}) lead to certain generalized Sinkhorn-type fixed point iterations \cite{karlsson2017generalized,cuturi2013sinkhorn} whose solutions yield the proximal updates $\left(\bm{\hat{\phi}}_k, \bm{\varpi}_{k}\right)$ for $k\in\mathbb{N}$. For algorithmic details, we refer the readers to \cite[Sec. III.B]{caluya2019gradient}. The convergence guarantees for such iterations follow from the contractive properties of cone preserving nonlinear maps that arise as first order conditions of optimality for (\ref{GradSinkhornProx}), and are detailed in \cite[Sec. III.C]{caluya2019gradient}.

For the mixed conservative-dissipative drift case, the \emph{uncontrolled} SDE associated with (\ref{transPDEdegen}) is
\begin{align}
\begin{pmatrix}
\differential\bm{\xi}\\
\differential\bm{\eta}	
\end{pmatrix} = \begin{pmatrix}
 \bm{\eta}\\
 -\nabla_{\bm{\xi}}V\left(\bm{\xi}\right) - \kappa\bm{\eta}	
 \end{pmatrix}\differential t + \sqrt{2\epsilon\kappa}\begin{pmatrix} 
 \bm{0}_{m\times m}\\
 \bm{I}_{m\times m}	
 \end{pmatrix}\differential\bm{w}.
\label{MixedUncontrolledSDE}	
\end{align} 
As before, for $k\in\mathbb{N}$, we denote the Euler-Maruyama updates associated with (\ref{MixedUncontrolledSDE}) with time step-sizes $\tau$ and $\sigma$, respectively, as $\{\bm{\xi}^{i}(t_{k-1}), \bm{\eta}^{i}(t_{k-1})\}_{i=1}^{N}\mapsto \{\bm{\xi}^{i}(t_{k}),\bm{\eta}^{i}(t_{k})\}_{i=1}^{N}$ and $\{\bm{\xi}^{i}(s_{k-1}), \bm{\vartheta}^{i}(s_{k-1})\}_{i=1}^{N}\allowbreak\mapsto \{\bm{\xi}^{i}(s_{k}),\bm{\vartheta}^{i}(s_{k})\}_{i=1}^{N}$. To reduce notational overload, let us recall the shorthands we used in (\ref{Wtilde})-(\ref{Defs}): $\bx=(\bm{\xi},\bm{\eta})^{\top}$, $\bm{y}=(\bm{\xi},\bm{\vartheta})^{\top}$, and write the Sinkhorn regularized proximal recursions for (\ref{MixedProx}) in vector form (see \cite[Sec. V.B]{caluya2019gradient}):
{\small{\begin{subequations}
	\begin{align}
&\bm{\hat{\phi}}_k =\: \underset{\bm{\hat{\phi}}} \argmin  \bigg\{ \underset{\bm{M} \in  \Pi(\bm{\hat{\phi}}_{k-1}, \bm{\hat{\phi}})}  \min \frac{1}{2} \bigg\langle\!\bm{S}_{\tau}\left(\{\bx^{i}(t_{k-1}),\bx^{i}(t_{k})\}_{i=1}^{N}\right),\bm{M}\!\bigg\rangle \nonumber\\
&+ \gamma \langle \bm{M}, \log \bm{M} \rangle +  \tau \bigg\langle\! \bigg\{\!\frac{1}{2}\|\bm{\eta}^{i}(s_{k-1})\|_{2}^{2}\bigg\}_{i=1}^{N}\!\!\!\!\!\! + \epsilon\log\bm{\varpi},\bm{\varpi} \!\bigg\rangle\! \bigg \},\label{MixedSinkhornProxphihat}\\
&\bm{\varpi}_k =\: \underset{\bm{\varpi}} \argmin   \bigg\{ \underset{\bm{M} \in  \Pi(\bm{\varpi}_{k-1}, \bm{\varpi})}  \min \frac{1}{2} \bigg\langle\! \bm{S}_{\sigma}\left(\{\bm{y}^{i}(s_{k-1}),\bm{y}^{i}(s_{k})\}_{i=1}^{N}\right),\bm{M}\!\bigg\rangle \nonumber\\
&+ \gamma \langle \bm{M}, \log \bm{M} \rangle + \sigma \bigg\langle\! \bigg\{\!\frac{1}{2}\|\bm{\vartheta}^{i}(s_{k-1})\|_{2}^{2}\bigg\}_{i=1}^{N}\!\!\!\!\!\!\! + \epsilon\log\bm{\varpi},\bm{\varpi} \!\bigg\rangle\! \bigg \},	\label{MixedSinkhornProxvarpi}
\end{align}
\label{MixedSinkhornProx}
\end{subequations}}}
wherein for $k\in\mathbb{N}$, the $(i,j)$\textsuperscript{th} component of the matrix $\bm{S}_{\tau}\in\mathbb{R}^{N\times N}$ in (\ref{MixedSinkhornProxphihat}) equals $s_{\tau}\left(\bm{x}^{i}(t_{k-1}),\bm{x}^{j}(t_{k})\right)$, and $s_{\tau}(\cdot,\cdot)$ is given by (\ref{Defstau}). Similarly, the $(i,j)$\textsuperscript{th} component of the matrix $\bm{S}_{\sigma}\in\mathbb{R}^{N\times N}$ in (\ref{MixedSinkhornProxvarpi}) equals $s_{\sigma}\left(\bm{y}^{i}(s_{k-1}),\bm{y}^{j}(s_{k})\right)$, and $s_{\sigma}(\cdot,\cdot)$ is given by (\ref{Defssigma}).

\subsection{Overall Algorithm}\label{SubsecOverallAlgo}
We now bring together the ideas from Sections \ref{SectionReformulation} and \ref{subsecWassProxSchrodinger}, and outline the overall algorithm to solve the SBPs with gradient and mixed-conservative dissipative drifts via scattered point cloud-based computation. We perform a recursion over the pair $(\hat{\varphi}_{0},\varphi_{1})$, and the computational steps for the same are:\\
 \noindent\textbf{Step 1.} Guess $\hat{\varphi}_{1}(\bx)$ (everywhere positive).\\
  \noindent\textbf{Step 2.} Compute $\varphi_{1}(\bx)=\rho_1(\bx)/\hat{\varphi}_{1}(\bx)$.\\
  \noindent\textbf{Step 3.} For the case of gradient drift, compute 
  \[p(\bx,s=0)\allowbreak =\allowbreak \varphi_{1}(\bx)\exp\left(-V(\bx)/\epsilon\right).\] For the case of mixed conservative-dissipative drift, compute 
  \[p(\bm{\xi},\bm{\vartheta},s=0) = \varphi_{1}(\bm{\xi},-\bm{\vartheta}) \exp\left( -\frac{1}{\epsilon} \left( \frac{1}{2} \lVert \bm{\vartheta}\rVert_{2}^{2} + V(\bm{\xi}) \right) \right).\]\\
  \noindent\textbf{Step 4.} Solve IVP (\ref{transPDE2}) or (\ref{transPDE2degen}) till $s=1$ to obtain $p(\bx,s=1)$.\\
    \noindent\textbf{Step 5.} For the case of gradient drift, compute
    \[\varphi_{0}(\bx) = p(\bx,s=1)\exp\left(V(\bx)/\epsilon\right).\]
    For the case of mixed conservative-dissipative drift, compute
    \[\varphi_{0}(\bm{\xi},\bm{\eta}) =p(\bm{\xi},\bm{\vartheta},s=1)\exp\left(\frac{1}{\epsilon}\left(\frac{1}{2}\|\bm{\eta}\|_{2}^{2} + V(\bm{\xi})\right)\right).\]\\
    \noindent\textbf{Step 6.} Compute $\hat{\varphi}_{0}(\bx) = \rho_{0}(\bx)/\varphi_{0}(\bx)$.\\
    \noindent\textbf{Step 7.} Solve IVP (\ref{transPDE1}) or (\ref{transPDE1degen}) till $t=1$ to obtain $\hat{\varphi}_{1}(\bx)$.\\
  \noindent\textbf{Step 8.} Repeat until the pair $(\hat{\varphi}_{0},\varphi_{1})$ has converged\footnote{We check whether the Wasserstein distances between the current and previous iterates for the pair $(\hat{\varphi}_{0},p(\bx,s=0))$, which is a proxy for the pair $(\hat{\varphi}_{0},\varphi_{1})$, are below some user-specified numerical tolerance. These Wasserstein distances are computed by solving linear programs, i.e.,  discrete versions of (41).}.\\

That such a recursion will converge, follows from the fact \cite[Proposition 1 in Sec. III]{chen2016entropic} that the recursion is in fact contractive in Hilbert's projective metric \cite{hilbert1895gerade,bushell1973hilbert} provided (i) the endpoint PDFs have compact supports, and (ii) the transition probability densities for (\ref{transPDE}) and (\ref{transPDEdegen}) are continuous and positive. Under the stated assumptions on the function $V$, the transition densities indeed satisfy these properties; see Appendix \ref{AppendinxRegularity}. Since our computational framework involves finite set of scattered points, the compactness condition also holds. Once the \emph{endpoint} Schr\"{o}dinger factors are found as outlined in Steps 1--8 above, we then compute the Schr\"{o}dinger factors at any time $t$, i.e., the pair $(\hat{\varphi}(\bx,t),\varphi(\bx,t))$ using the IVPs (\ref{transPDE}) or (\ref{transPDEdegen}).

\begin{algorithm}[t]
    \caption{\textproc{ComputeFactorsSBP} using \textproc{ProxRecur} to compute the Schr\"{o}dinger factors $(\bm{\hat{\phi}}^{\rm{transient}}_k,\bm{\phi}^{\rm{transient}}_k)$ for  gradient drift}
    \label{AlgoComputeFactorsSBP}
    \begin{algorithmic}[1] 
        \Procedure{ComputeFactorsSBP}{$\bm{\rho}_0$, $\bm{\rho}_1$, $\gamma$, $\epsilon$, $\tau$, $\sigma$, $N$, $V(\cdot)$, $\rm{numSteps}$, $\rm{tol}_{\rm{SB}}$, $\rm{maxIter}_{\rm{SB}}$, $\rm{tol}_{\rm{PR}}$, $\rm{maxIter}_{\rm{PR}}$}
    \State $\bm{\hat{\phi}}_1\gets \left[ {\rm{rand}}_{N \times 1},\bm{0}_{N\times( \rm{maxIter}_{\rm{SB}}-1)}\right]$ \Comment{Step 1}
            \State $\bm{\phi}_0 \gets  \left[ \bm{0}_{N \times \rm{maxIter}_{\rm{SB}} }\right]$ \Comment{initialize}
            \State $\bm{\hat{\phi}}_0 \gets \left[ \bm{0}_{N \times \rm{maxIter}_{\rm{SB}} }\right]$
            \State $\bm{\phi}_1 \gets \left[ \bm{0}_{N \times \rm{maxIter}_{\rm{SB}} }\right] $
             \State $\bm{p}_0 \gets \left[ \bm{0}_{N \times \rm{maxIter}_{\rm{SB}} }\right] $
              \State $\bm{p}_1 \gets \left[ \bm{0}_{N \times \rm{maxIter}_{\rm{SB}} }\right] $
              \State $\bm{p}^{\rm{temp}} \gets \left[ \bm{0}_{N\times \rm{numSteps}+1} \right]$
         	 \State $\bm{\hat{\phi}}^{\rm{temp}} \gets \left[ \bm{0}_{N\times \rm{numSteps}+1} \right]$
           
            \State $\ell=1$ \Comment{iteration index for \textproc{ComputeFactorsSBP}}
            \While{$\ell \leq \rm{maxIter}_{\rm{SB}}$}
                \State $\bm{\phi}_1(:,\ell+1) = \bm{\rho}_1\oslash \bm{\hat{\phi}}_1(:,\ell)$ \Comment{Step 2}
                \State $\bm{p}_0(:,\ell+1)= \bm{\phi}_1(:,\ell+1) \odot \exp(-V(\{\bx^{i}\}_{i=1}^{N})/\epsilon) $ \Comment{\parbox[t]{.3\linewidth}{Step 3}}
                 \State $\bm{p}^{\rm{temp}}(:,1) = \bm{p}_0(:,\ell+1)$
                \For{$i \gets 1$ to $\rm{numSteps}$}
                 \State $\bm{p}^{\rm{temp}}(:,i+1) \gets \textproc{ProxRecur}(\bm{p}^{\rm{temp}}(:,i),\gamma,$
                 \StatexIndent $\epsilon, \sigma, N,V(\cdot),\rm{tol}_{\rm{PR}}$,$\rm{maxIter}_{\rm{PR}})$
                 \EndFor
                 \State $\bm{p}_1(:,\ell+1) \gets \bm{p}^{\rm{temp}}(:,\rm{numSteps}+1)$ \Comment{Step 4}
                 \State $\bm{\phi}_0(:,\ell+1) \gets \bm{p}_1(:,\ell+1) \odot \exp(V\left(\{\bx^{i}\}_{i=1}^{N}\right)/\epsilon) $ \Comment{Step 5}
                \State $\bm{\hat{\phi}}_0(:,\ell+1) \gets \bm{\rho}_0 \oslash \bm{\phi}_0(:,\ell+1)$ \Comment{Step 6}
                \State $\bm{\hat{\phi}}_{\rm{temp}}(:,1)\gets \bm{\hat{\phi}}_0(:,\ell+1)$
                \For{$j \gets 1$ to $\rm{numSteps} $}
                \State $\bm{\hat{\phi}}^{\rm{temp}}(:,j+1) \gets \textproc{ProxRecur}(\bm{\hat{\phi}}^{\rm{temp}}(:,j),\gamma,$
                \StatexIndent $\epsilon,\tau,N,V(\cdot),\rm{tol}_{\rm{PR}},\rm{maxIter}_{\rm{PR}})$
                \EndFor
                   \State $\bm{\hat{\phi}}_1(:,\ell+1) \gets \bm{\hat{\phi}}^{\rm{temp}}(:,\rm{numSteps}+1)$ \Comment{Step 7}
                \StatexIndent
                \If{$W^2( \bm{p}_0(:,\ell+1),\bm{p}_0(:,\ell) ) <  {\rm{tol}}_{{\rm{SB}}} \And W^2(\bm{\hat{\phi}}_0(:,\ell+1),\bm{\hat{\phi}}_0(:,\ell)) < \rm{tol}_{\rm{SB}}$ } 
                \State break
                \Else
                \State $\ell \gets \ell + 1$ \Comment{Step 8}
                \EndIf
            \EndWhile\label{euclidendwhile}
            \State \textbf{return} $\bm{\hat{\phi}}_0(:,\ell),\bm{\phi}_1(:,\ell)$\Comment{converged endpoint pair}
            \StatexIndent
            \State $\bm{\hat{\phi}}^{\rm{transient}}_{k=1} \gets \bm{\hat{\phi}}_{0}(:,\ell)$ \Comment{initialize}
            \State $\bm{p}^{\rm{transient}}_{k=1}\gets \bm{\phi}_{1}(:,\ell) \odot \exp(-V\left(\{\bx^{i}\}_{i=1}^{N}\right)_1/\epsilon)$ \Comment{initialize}
             \For{$k \gets 1$ to $\rm{numSteps}$}
           \State $ \bm{\hat{\phi}}^{\rm{transient}}_{k+1}\gets \textproc{ProxRecur} (\bm{\hat{\phi}}^{\rm{transient}}_{k},\gamma,\epsilon, \tau, N,V(\cdot),$
           \StatexIndent $\rm{tol}_{\rm{PR}},\rm{maxIter}_{\rm{PR}})$
            \State $\bm{p}^{\rm{transient}}_{k+1} \gets \textproc{ProxRecur}(\bm{p}^{\rm{transient}}_{k},\gamma,\epsilon, \sigma,N,V(\cdot),$
            \StatexIndent $\rm{tol}_{\rm{PR}},\rm{maxIter}_{\rm{PR}})$
            \State $\bm{\phi}^{\rm{transient}}_{k+1}\gets\bm{p}^{\rm{transient}}_{k+1}\odot \exp(V\left(\{\bx^{i}\}_{i=1}^{N}\right)/\epsilon)$
           \EndFor
        \State\textbf{return} $(\bm{\hat{\phi}}^{\rm{transient}}_k,\bm{\phi}^{\rm{transient}}_k)$ \!\!\Comment{transient Schr\"{o}dinger factors}   
        \EndProcedure
    \end{algorithmic}
\end{algorithm}

For scattered point cloud-based computation with sample size $N$, the aforesaid steps lead to a recursion over the pair of vectors $(\bm{\hat{\phi}}_0,\bm{\phi}_1)$, each of these vectors being of size $N\times 1$. We refer to this proposed procedure as \textproc{ComputeFactorsSBP} which takes the $N\times 1$ endpoint PDF vectors $\bm{\rho}_0,\bm{\rho}_1$ (i.e., the endpoint joint PDF values evaluated at the scattered sample state vectors) as input, and returns the converged pair $(\bm{\hat{\phi}}_0,\bm{\phi}_1)$. As detailed in Algorithm 1, \textproc{ComputeFactorsSBP} performs the discrete versions of the Steps 1--8 above, followed by computation of the transient Schr\"{o}dinger factors. For conceptual clarity,  in Algorithm 1, we list the steps for the case of gradient drift; adapting Algorithm 1 to the case of mixed conservative-dissipative drift is straightforward following the modifications in Steps 3, 4, 5 and 7 mentioned before. 

The procedure \textproc{ComputeFactorsSBP} involves a main iteration where we guess an initial vector $\bm{\hat{\phi}}_1$, and repeatedly update the vectors $\bm{\hat{\phi}}_0,\bm{\hat{\phi}}_1,\bm{\phi}_0,\bm{\phi}_1$ using a secondary function called \textproc{ProxRecur} that performs the proximal recursions (\ref{GradSinkhornProx}) or (\ref{MixedSinkhornProx}) for solving the IVPs (\ref{transPDE}) or (\ref{transPDEdegen}) in Steps 4 and 7 above. The main iteration requires the parameters $\gamma,\epsilon,\tau,\sigma,N$, the function $V(\cdot)$ described before, and $\rm{numSteps}$ which is the number of time-steps for the proximal update. Furthermore, we use the two pairs of parameters 
$(\rm{tol}_{\rm{SB}},\rm{maxIter}_{\rm{SB}}),(\rm{tol}_{\rm{PR}} ,\rm{maxIter}_{\rm{PR}})$: the first pair specifies the numerical tolerance and the maximum number of iterations for the main iteration in \textproc{ComputeFactorsSBP}; the second pair specifies the same for the Sinkhorn algorithm used in \textproc{ProxRecur}. Notice that the procedure  \textproc{ProxRecur} subsumes the Euler-Maruyama scheme to update the states $\{\bx^{i}\}_{i=1}^{N}$; see \cite[Sec. III.B]{caluya2019gradient}.

Once the converged vector pair $(\bm{\hat{\phi}}_0,\bm{\phi}_1)$ is obtained, we again invoke \textproc{ProxRecur} to solve the IVP pair (\ref{transPDE}) or (\ref{transPDEdegen}) via (\ref{GradSinkhornProx}) or (\ref{MixedSinkhornProx}), and compute the discrete transient solutions ($\bm{\hat{\phi}}_{k},\bm{\varpi}_{k}$). We remind the readers that since the function $p$ is integrated in the time coordinate $s$, its discrete version $\bm{\varpi}_{k}$ is really evaluated at $t= 1- k\tau$, i.e., $\bm{\varpi}_{k} \approx p(\bx,t=1-k\tau)$ for $k\in\mathbb{N}$. We next employ the discrete versions of the mappings $p \mapsto \varphi$, given by
\begin{align*}
\varphi(\bx,t) &= p(\bx,s)\exp\left(V(\bx)/\epsilon\right),\\
\varphi(\bm{\xi},\bm{\eta},t) &= p(\bm{\xi},\bm{\vartheta},s)\exp\left(\frac{1}{\epsilon}\left(\frac{1}{2}\|\bm{\eta}\|_{2}^{2} + V(\bm{\xi})\right)\right),	
\end{align*}
to compute $\bm{\varpi}_{k}\mapsto \bm{\phi}_{k}$, i.e., $\left(\bm{\hat{\phi}}_{k},\bm{\varpi}_{k}\right) \mapsto \left(\bm{\hat{\phi}}_{k},\bm{\phi}_{k}\right)$. This completes the computation of the transient Schr\"{o}dinger factors $\left(\bm{\hat{\phi}}_{k},\bm{\phi}_{k}\right)$ as a pair of weighted scattered point clouds, i.e., as function values $(\hat{\varphi}(\bx,t),\varphi(\bx,t))$ evaluated at the state coordinates that are updated by the Euler-Maruyama scheme.

To compute the discrete optimal pair $(\bm{\rho}^{\rm{opt}}_{k},\bm{u}^{\ropt}_{k})$ from $\left(\bm{\hat{\phi}}_{k},\bm{\phi}_{k}\right)$, we need to evaluate the Hadamard product
\begin{equation} 
\bm{\rho}^{\rm{opt}}_{k} = \bm{\phi}_{{\rm{numSteps}}+1-k} \odot \bm{\hat{\phi}}_{k},
\label{HadadmardProductDiscreteRhoOpt}
\end{equation} 
which in our case, is not well-defined as is, since $\bm{\phi}^{\rm{trans}}_{{\rm{numSteps}}+1-k}$ and $\bm{\hat{\phi}}^{\rm{trans}}_{k}$ are evaluated on different state  coordinates (i.e., are supported on different finite sets). On the other hand, the optimal feedback control
\begin{equation} \label{Discretecontrol}
\bm{u}^{\rm{opt}}_{k} = 2 \epsilon \bm{B}^{\top}\nabla \log\bm{\phi}_{k},
\end{equation}
requires computing a gradient w.r.t. the state, which again is not well-defined for scattered data. To circumvent these issues, we perform scattered data interpolation for the transient Schr\"{o}dinger factors $\left(\bm{\hat{\phi}}_{k},\bm{\phi}_{k}\right)$ via \emph{multiquadric} scheme \cite{hardy1971multiquadric,hardy1975research} that uses radial basis functions to fit a surface to the weighted scattered point cloud. The interpolation requires the vectors $(\bm{\phi}_k,\bm{\hat{\phi}}_k)$, their corresponding state space coordinates, a set of user-provided query points, and returns the interpolated values evaluated at the query points. Using the same set of query points to interpolate the two point clouds for the pair $(\bm{\phi}_k,\bm{\hat{\phi}}_k)$, we perform the element-wise multiplication (\ref{HadadmardProductDiscreteRhoOpt}). As for (\ref{Discretecontrol}), we use standard finite difference techniques to compute the gradient of the logarithm of interpolated values of $\bm{\phi}_k$. Alternatively, (\ref{Discretecontrol}) could be evaluated via kernel density estimation.


\section{Numerical Examples}\label{SectionNumericalExamples}


\subsection{SBP Example with Gradient Drift}\label{GradSBPexample}
We consider an instance of the SBP given in Section \ref{SubsectionSchrodingerSystemwithGradientDrift} with $\bx\in \mathbb{R}^2$, and $V(x_1,x_2) = \frac{1}{4}(1+x_1^4) + \frac{1}{2}(x_2^2-x_1^2)$. The controlled prior dynamics is
\begin{align} 
\begin{pmatrix} \differential x_1\\ \differential x_2 \end{pmatrix}  = -\nabla V(x_1,x_2) \:\differential t  + \begin{pmatrix} u_1\\  u_2\end{pmatrix}  \:\differential t +\sqrt{2\epsilon}\begin{pmatrix} \differential w_1 \\ \differential w_2 \end{pmatrix}. 
\label{SDEexample2dcont}
\end{align}
The control objective is to steer the prescribed joint PDF of the initial condition $\bx(t=0) \sim \rho_0 = \mathcal{N}(\bm{\mu}_0,\bm{\Sigma}_0)$ to the prescribed joint PDF of the terminal condition $\bx(t=1) \sim \rho_1 = c_{1} \mathcal{N}(\bm{\mu}_{1},\bm{\Sigma}_{1})  + c_{2} \mathcal{N}(\bm{\mu}_{2},\bm{\Sigma}_{2})$ over $t\in[0,1]$, subject to (\ref{SDEexample2dcont}), while minimizing the control effort (\ref{SBwGradDriftObj}). Here, we fix
\begin{align*}
&\bm{\mu}_0 = (-2,0)^{\!\!\top}, \: \bm{\Sigma}_0 =  {\rm{diag}}(0.8,0.7), \: c_{1}=c_{2}=0.5,\:\bm{\mu}_{1} = (1.5,2)^{\!\!\top},\\
&\bm{\mu}_{2} = (1.5,-2)^{\!\!\top},\:\bm{\Sigma}_{1} =  {\rm{diag}}(0.5,0.8),\:\bm{\Sigma}_{2} =  {\rm{diag}}(0.7,0.8).
\end{align*}
Notice that in the absence of control ($\bm{u}\equiv0$), the transient joint PDFs of (\ref{SDEexample2dcont})
tend to the stationary solution $\rho_{\infty}\propto\exp\left(-V(x_1,x_2)/\epsilon\right)$ which has two modes along the horizontal axis; see \cite[Fig. 8]{caluya2019gradient}. In contrast, the prescribed terminal bimodal PDF $\rho_{1}$, specified as a mixture of Gaussians, has two modes along the vertical axis.

\begin{figure*}[t]
\centering
\includegraphics[width=.87\linewidth]{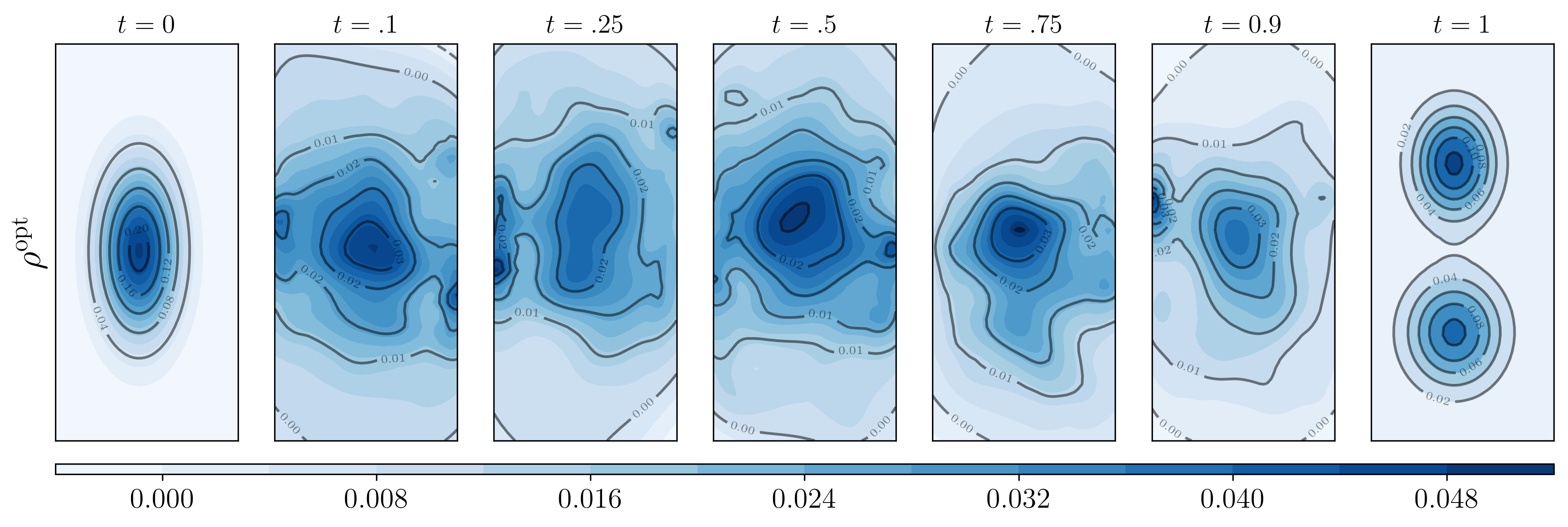}
\vspace*{-0.10in}
\caption{\small{For the SBP in Section \ref{GradSBPexample}, shown here are the contour plots of the optimal controlled transient joint state PDFs $\rho^{{\rm{opt}}}(\bx,t)$, $t\in[0,1]$, along with the endpoint joint PDFs $\rho_0(\bx),\rho_1(\bx)$. Each subplot corresponds to a different snapshot in time; all subplots are plotted on the domain $[-4,4]\times[-6,6]$. The color denotes the joint PDF value; see colorbar (dark hue = high, light hue  = low).}}
\label{FigGradSBPOptRho}
\vspace*{-0.1in}
\end{figure*}

\begin{figure*}[t]
\centering
\includegraphics[width=.87\linewidth]{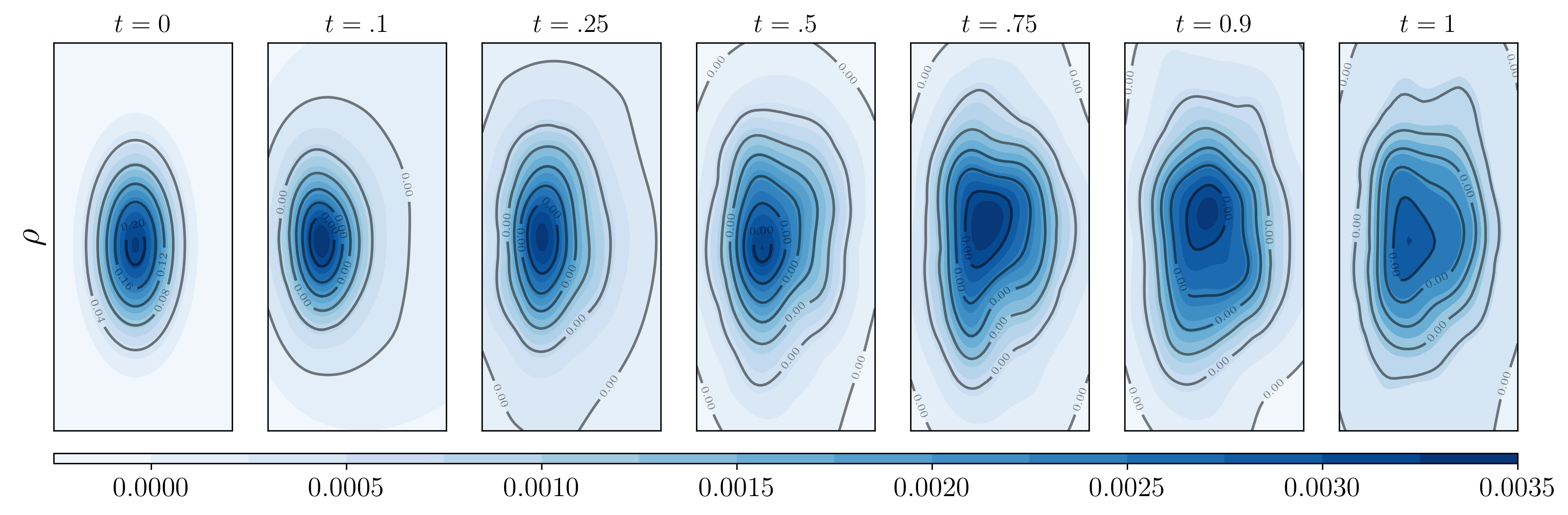}
\vspace*{-0.10in}
\caption{\small{The contour plots of the uncontrolled  ($\bm{u}\equiv0$) transient joint state PDFs $\rho(\bx,t)$, $t\in[0,1]$, for (\ref{SDEexample2dcont}) starting from the initial joint state PDF $\rho_0$ given in Section \ref{GradSBPexample}. Each subplot corresponds to a different snapshot in time. Each subplot corresponds to a different snapshot in time; all subplots are plotted on the domain $[-4,4]\times[-6,6]$. The color denotes the joint PDF value; see colorbar (dark hue = high, light hue  = low).}}
\label{FigGradSBPUncontrolledRho}
\vspace*{-0.1in}
\end{figure*}

\begin{figure*}[t]
\centering
\includegraphics[width=.87\linewidth]{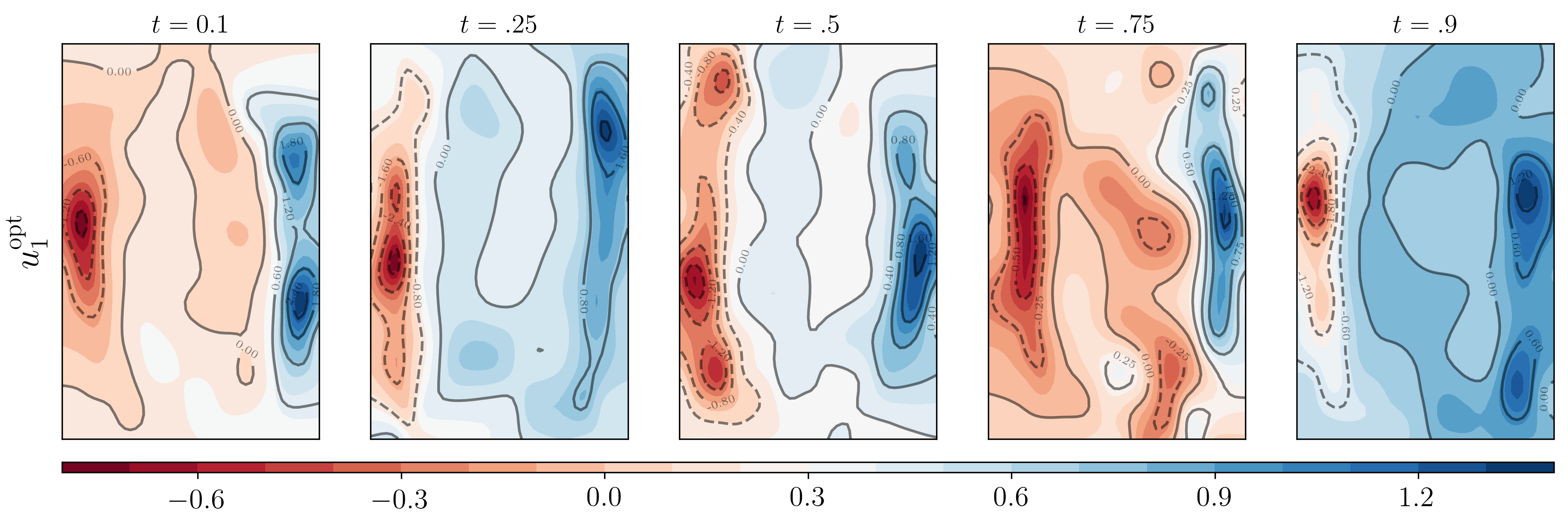}
\vspace*{-0.10in}
\caption{\small{For the SBP in Section \ref{GradSBPexample}, shown here are the contour plots of $u_1^{\ropt}(\bx,t)$, the first component of the optimal feedback control.  Each subplot is plotted on the domain $[-4,4]\times[-6,6]$. The color (blue = high, red = low) denotes the value of  $u_1^{\ropt}$ at each snapshot in time.}}
\label{FigGradSBPOptU1}
\vspace*{-0.15in}
\end{figure*}

\begin{figure*}[t]
\centering
\includegraphics[width=.87\linewidth]{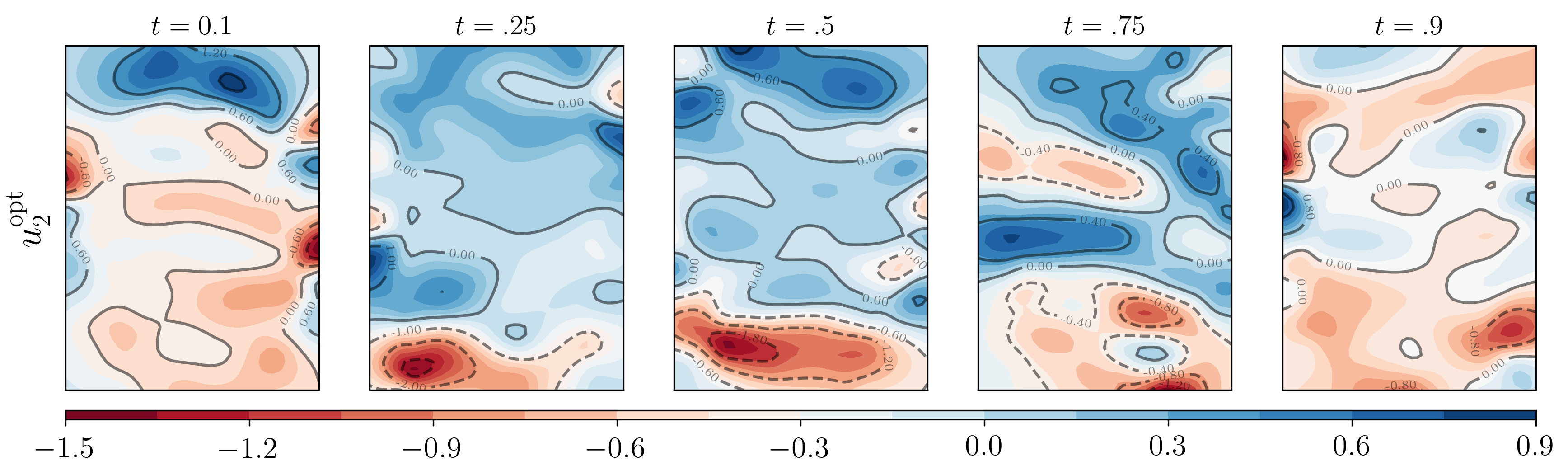}
\vspace*{-0.10in}
\caption{\small{For the SBP in Section \ref{GradSBPexample}, shown here are the contour plots of $u_2^{\ropt}(\bx,t)$, the second component of the optimal feedback control.  Each subplot is plotted on the domain $[-4,4]\times[-6,6]$. The color (blue = high, red = low) denotes the value of  $u_2^{\ropt}$ at each snapshot in time.}}
\label{FigGradSBPOptU2}
\vspace*{-0.15in}
\end{figure*}

\begin{figure*}[t]
\centering
\includegraphics[width=.87\linewidth]{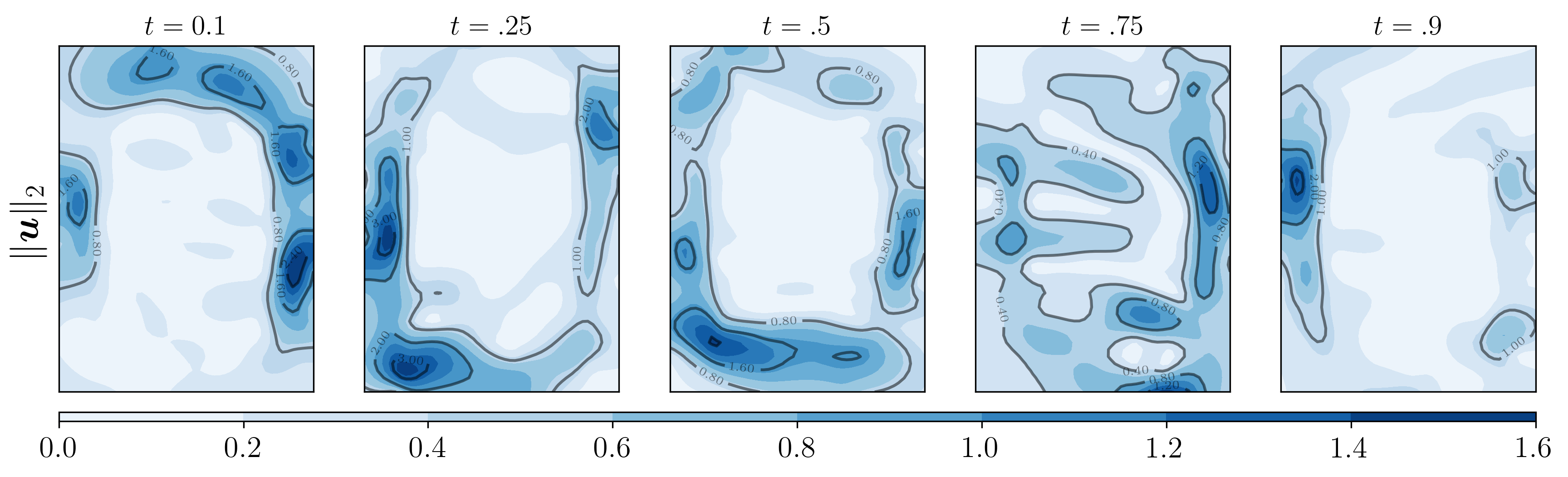}
\vspace*{-0.10in}
\caption{\small{For the SBP in Section \ref{GradSBPexample}, shown here are the contour plots for the magnitude  (dark hue = high, light hue  = low; see colorbar) of the optimal feedback control, i.e., $\lVert \bm{u} \rVert_2= \sqrt{(u_1^{\ropt}(\bx,t))^2 +(u_2^{\ropt}(\bx,t))^2 }$.}}
\label{FigGradSBPOptUmag}
\vspace*{-0.15in}
\end{figure*}

Fig. \ref{FigGradSBPOptRho} shows the evolution of the optimal controlled joint PDF $\rho^{{\rm{opt}}}(\bx,t)$ obtained by solving the SBP using Algorithm 1 given in Section \ref{SubsecOverallAlgo}. For numerical simulation, we set $\epsilon = 6$, $\gamma = 0.5$, $\tau=\sigma=10^{-3}$, $N= 500$, $\rm{numSteps} = 1000$, $\rm{tol}_{\rm{SB}}=0.1$, $\rm{maxIter}_{\rm{SB}}=\rm{maxIter}_{\rm{PR}}=500$, and $\rm{tol}_{\rm{PR}}=10^{-3}$. In Fig. \ref{FigGradSBPOptRho}, each subplot corresponds to a different snapshot in time; each subplot is plotted over the domain $[-4,4]\times[-6,6]$. For the purpose of comparison, Fig. \ref{FigGradSBPUncontrolledRho} shows the uncontrolled joint state PDF evolution starting from the same initial PDF $\rho_{0}$. The components of the optimal feedback control $\bm{u}^{{\rm{opt}}}(\bx,t)$ obtained from the SBP solution, are shown in Figs. \ref{FigGradSBPOptU1} and \ref{FigGradSBPOptU2}. Fig. \ref{FigGradSBPOptUmag} depicts the magnitude of $\bm{u}^{\ropt}$.

\subsection{SBP Example with Mixed Conservative-Dissipative Drift}\label{MixedSBPexample}
We next consider an instance of the SBP given in Section \ref{SubsectionSchrodingerSystemwithMixedDrift} with $n=2m=2$, i.e., $\bx = (\xi,\eta)^{\top}\in\mathbb{R}^{2}$, and $V(\xi) = 5\xi^{4}$. In other words, the controlled prior dynamics is
\begin{align} 
\begin{pmatrix} \differential \xi\\ \differential \eta \end{pmatrix}  = \begin{pmatrix} \eta\\  -\frac{\partial}{\partial\xi}5\xi^{4}-\kappa\eta + u\end{pmatrix}  \:\differential t +\sqrt{2\epsilon\kappa}\begin{pmatrix} 0 \\ \differential w \end{pmatrix}. 
\label{SDEexampleMixed}
\end{align}
Notice that the function $V$ satisfies the conditions mentioned in Section \ref{SubsectionSchrodingerSystemwithMixedDrift}. We use the same endpoint PDFs $\rho_{0},\rho_{1}$ as in Section \ref{GradSBPexample}, and solve the SBP using Algorithm 1 with the Steps 3, 4, 5 and 7 modified for the mixed conservative-dissipative case, as mentioned before. Notice that in the absence of control  ($u\equiv0$), the transient joint PDFs of (\ref{SDEexampleMixed})
tend to the stationary solution $\rho_{\infty}\propto\exp\left(-\frac{1}{\epsilon}\left(\frac{1}{2}\eta^{2} + V(\xi)\right)\right)$.

\begin{figure*}[t]
\centering
\includegraphics[width=.87\linewidth]{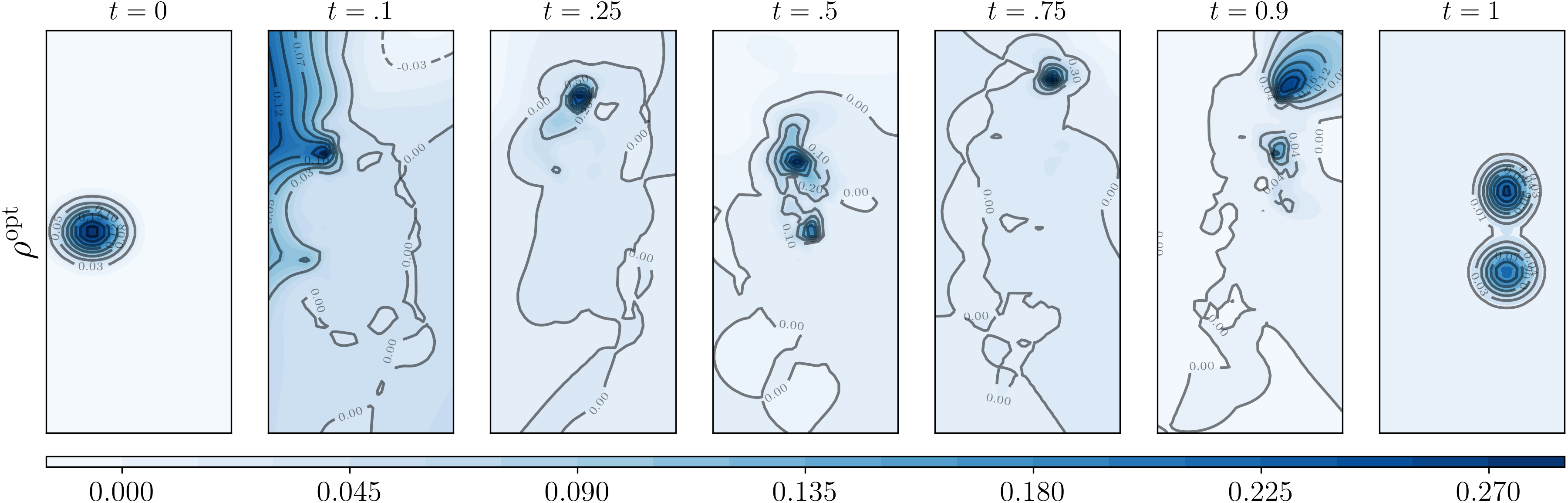}
\caption{\small{For the SBP in Section \ref{MixedSBPexample}, shown here are the contour plots of the optimal controlled transient joint state PDFs $\rho^{{\rm{opt}}}(\bx,t)$, $t\in[0,1]$, along with the endpoint joint PDFs $\rho_0(\bx),\rho_1(\bx)$. Each subplot corresponds to a different snapshot in time; all subplots are plotted on the domain $[-4,4]\times[-10,10]$. The color denotes the joint PDF value; see colorbar (dark hue = high, light hue  = low).}}
\label{rhotransdegen}
\vspace*{-0.15in}
\end{figure*}

\begin{figure*}[t]
\centering
\includegraphics[width=.87\linewidth]{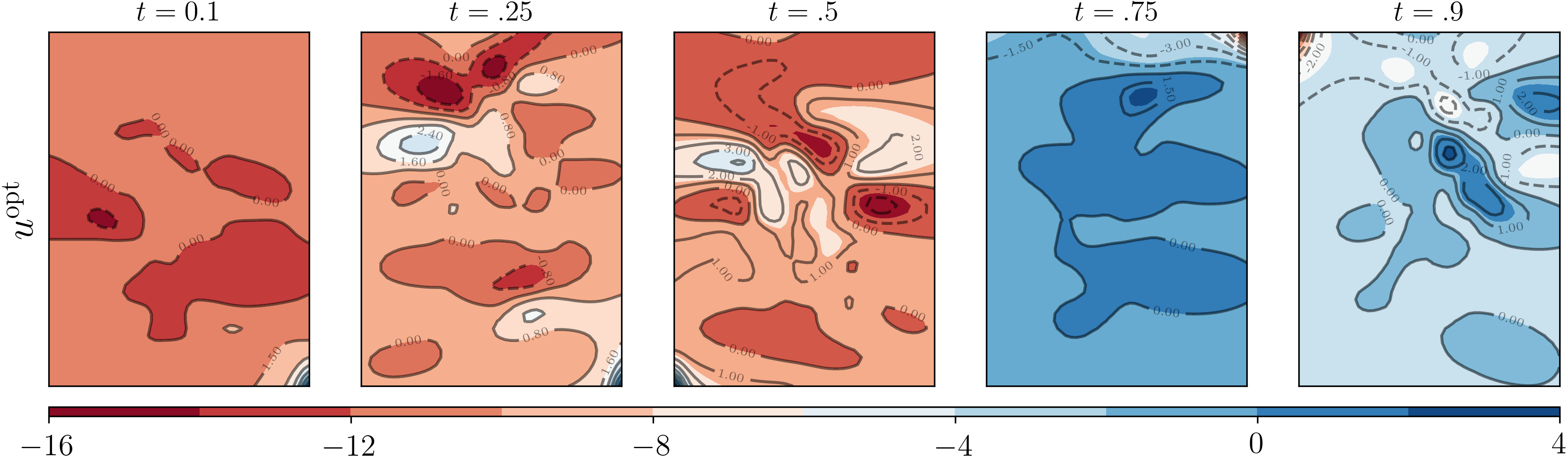}
\caption{\small{For the SBP in Section \ref{MixedSBPexample}, shown here are the contour plots of the optimal feedback control $u^{\ropt}(\bx,t)$.  Each subplot is plotted on the domain $[-4,4]\times[-10,10]$. The color (blue = high, red = low) denotes the value of  $u^{\ropt}$ at each snapshot in time; see colorbar.}}
\label{uoptdegen}
\vspace*{-0.15in}
\end{figure*}

Fig. \ref{rhotransdegen} shows the evolution of the optimal controlled joint state PDF $\rho^{{\rm{opt}}}(\bx,t)$ obtained from the SBP solution, where each subplot is plotted over the domain $[-4,4]\times[-10,10]$. The parameters used in our numerical simulation are: $\epsilon = 5$, $\kappa = 0.5$, $\gamma = 0.5$, $\tau=\sigma = 10^{-3}$, $N = 100$, $\rm{numSteps} = 1000$, $\rm{tol}_{\rm{SB}}=0.1$, $\rm{maxIter}_{\rm{SB}}=\rm{maxIter}_{\rm{PR}}=500$, and $\rm{tol}_{\rm{PR}}=10^{-3}$. The optimal feedback control $\bm{u}^{{\rm{opt}}}(\bx,t)$ obtained from our SBP solution is shown in Fig. \ref{uoptdegen}.

\section{Conclusions}\label{SectionConclusions}
We address the problem of minimum energy finite horizon steering of state density between two prescribed endpoint densities via feedback control, subject to trajectory-level dynamics with nonlinear drift. This is a generalization of the classical SBP formulated by Schr\"{o}dinger in the 1930s, and has many potential applications such as shaping of biological and robotic swarms through feedback synthesis. We derive optimality conditions for the case of a generic nonlinear drift, and show that two specific cases, viz. gradient and mixed-conservative dissipative drifts, are particularly amenable for computational purpose, based on certain infinite-dimensional proximal recursions that exploit the structure of these nonlinearities. These ideas make contact with the theory of Wasserstein gradient flow, and optimal mass transport. Building on our previous work \cite{caluya2019gradient}, we design proximal algorithms for solving the density steering problems in these two cases of practical interest. Numerical examples are provided to illustrate the proposed framework. Our main contribution is to make algorithmic advances in solving the variants of SBPs that are of interest to the systems-control community. From a probabilistic perspective, while there exists a substantial body of work \cite{hotz1987covariance,halder2016finite,bakolas2018finite,okamoto2018optimal} in steering the second order state statistics (also known as ``covariance control"), steering the joint state PDF subject to controlled nonlinear dynamics is a relatively new direction of research. The results of this article are expected to further the theoretical and algorithmic development for the same.

\appendix
\vspace*{-0.1in}
\subsection{Proof of Proposition \ref{propSBPnonlineargeneral}}\label{AppendixCondOpt}
\vspace*{-0.1in}
\noindent We start by rewriting the Lagrangian (\ref{SBwDriftLagrangian}) as 
{\small{\begin{equation}
\begin{aligned} \label{SBwDriftLagrangian2}
\mathcal{L}(\rho,\bm{u},\psi) &= \int_{0}^{1} \int_{\bR^n} \frac{1}{2} \lVert \bm{u}(\bx,t) \rVert_{2}^{2}\:\rho(\bx,t)  \: \differential\bx\:\differential t \\
&+  \underbrace{\int_{0}^{1} \int_{\bR^n}  \psi(\bx,t) \frac{\partial \rho}{\partial t}   \: \differential\bx\:\differential t }_{\text{term 1}}\\ 
&\hspace*{-0.6in}+  \!\!\underbrace{\int_{0}^{1}\!\!\int_{\bR^n}\!\!\psi(\bx,t)  \bigg\{\nabla \!\cdot\!\left(\rho\left(\bm{f}+\bm{B}(t)\bm{u}\right)\right) -\epsilon\langle\bm{D}(t),\hess\left(\rho\right)\rangle\bigg\}\: \differential\bx\:\differential t}_{\text{term 2}}.
\end{aligned}
\end{equation}}}
Next, we switch the order of integration, perform integration by parts w.r.t. $t$ in term 1, and integration by parts w.r.t. $\bx$ in term 2. Assuming the limits for $\parallel\bx\parallel_{2} \rightarrow \infty$ are zero, the Lagrangian (\ref{SBwDriftLagrangian2}) then simplifies to
\begin{align}
\mathcal{L}(\rho,\bm{u},\psi) & =  \int_{0}^{1} \int_{\bR^n}  \bigg( \frac{1}{2}\lVert \bm{u}(\bx,t) \rVert_{2}^{2}- \frac{ \partial \psi }{\partial t} \nonumber\\
 &\hspace*{-0.5in}-\langle\nabla \psi,\bm{f}+\bm{B}(t)\bm{u}\rangle - \epsilon \langle\bm{D}(t),\hess\left(\psi\right)\rangle \bigg) \rho(\bx,t)
  \: \differential\bx\:\differential t,
  \label{SBwDriftLagrangian3} 
\end{align}
where we have used the (two-fold) integration by parts:
\begin{align*}
&\displaystyle\int_{\mathbb{R}^{n}}\langle\bm{D}(t),\hess\left(\rho\right)\rangle\:\psi(\bx,t)\:\differential\bx\\
=& \displaystyle\sum_{i,j=1}^{n} \displaystyle\int_{\mathbb{R}^{n}} \displaystyle\frac{\partial^{2}\left(\bm{D}(t)\rho(\bx,t)\right)_{ij}}{\partial x_{i}\partial x_{j}}\:\psi(\bx,t)\:\differential\bx\\
=& -\displaystyle\sum_{i,j=1}^{n} \displaystyle\int_{\mathbb{R}^{n}}\dfrac{\partial\left(\bm{D}(t)\rho\right)_{ij}}{\partial x_{j}}\dfrac{\partial\psi}{\partial x_{i}}\:\differential\bx\\
=& \displaystyle\sum_{i,j=1}^{n}\displaystyle\int_{\mathbb{R}^{n}}\left(\bm{D}(t)\rho\right)_{ij}\dfrac{\partial^{2}\psi}{\partial x_{j}\partial x_{i}}\:\differential\bx\\
=& \displaystyle\int_{\mathbb{R}^{n}}\displaystyle\sum_{i,j=1}^{n}\left(\bm{D}(t)\rho\right)_{ij}\dfrac{\partial^{2}\psi}{\partial x_{j}\partial x_{i}}\:\differential\bx\\
=& \displaystyle\int_{\mathbb{R}^{n}} \langle\bm{D}(t),\hess\left(\psi\right)\rangle\rho(\bx,t)\: \differential\bx.
\end{align*}

\noindent Pointwise minimization of (\ref{SBwDriftLagrangian3}) w.r.t. $\bm{u}$ for a fixed PDF $\rho$, gives 
\begin{equation}
\bm{u}^{{\rm{opt}}}(\bx,t) = \bm{B}(t)^{\top}\nabla \psi(\bx,t).
\label{Provinguoptequalsgradpsi}   
\end{equation}
Substituting (\ref{Provinguoptequalsgradpsi}) back into (\ref{SBwDriftLagrangian3}), and equating the resulting expression to zero, we get the dynamic programming equation 
\begin{equation}\label{DPeqn}
\begin{aligned}
 \int_{0}^{1}\int_{\bR^{n}} \left( -\frac{\partial\psi}{\partial t} -\frac{1}{2} \lVert \nabla \bm{B}(t)^{\top}\psi \rVert_{2}^{2} - \langle\nabla\psi,\bm{f}\rangle \right. \\
 - \left.\right.\epsilon \langle\bm{D}(t),\hess\left(\psi\right)\rangle \bigg) \rho(\bx,t) \: \differential\bx \: \differential t = 0.
\end{aligned}
\end{equation}
For (\ref{DPeqn}) to hold for arbitrary $\rho$, we must have
\begin{equation}
\begin{aligned}
\frac{\partial \psi }{\partial t} + \frac{1}{2}\lVert \bm{B}(t)^{\top}\nabla\psi \rVert_{2}^{2}  + \langle\nabla\psi, \bm{f}\rangle = - \epsilon \langle\bm{D}(t),\hess\left(\psi\right)\rangle,
\end{aligned}
\end{equation}
which is the HJB PDE (\ref{HJBgenNonlinear}). The associated FPK PDE (\ref{FPKgenNonlinear}) results from substituting (\ref{Provinguoptequalsgradpsi}) into (\ref{SBFluidsFP}). The boundary conditions (\ref{SBwDriftBC}) follows from (\ref{SBFluidsBC}). This completes the proof.
\hfill\qedsymbol

\vspace*{-0.1in}
\subsection{Proof of Theorem \ref{ThmHopfCole}}\label{AppendixHopfColeProof}
\vspace*{-0.1in}
\noindent In (\ref{HopfColevarphi}), taking the gradient of $\varphi$ w.r.t. $\bx$, we get   
\begin{align}
\nabla \varphi = \frac{1}{2\epsilon}\exp\left(\frac{\psi}{2\epsilon}\right) \nabla \psi.
\label{nablaphiANDnablapsi}
\end{align}
Furthermore, 
\begin{align} 
&-\epsilon\langle\bm{D}(t),\hess\left(\varphi\right)\rangle = -\epsilon\displaystyle\sum_{i,j=1}^{n}(\bm{D}(t))_{ij}\dfrac{\partial^{2}}{\partial x_{i}\partial x_{j}} \exp\left(\psi/2\epsilon\right)\nonumber\\
 &= -\epsilon\dfrac{\exp\left(\psi/2\epsilon\right)}{2\epsilon}\bigg\{\! \displaystyle\sum_{i,j=1}^{n}(\bm{D}(t))_{ij}\left(\dfrac{\partial^{2}\psi}{\partial x_{i}\partial x_{j}} + \frac{1}{2\epsilon}\dfrac{\partial\psi}{\partial x_{i}}\dfrac{\partial\psi}{\partial x_{j}}\right) \!\!\bigg\}\nonumber\\
 &= \frac{1}{2\epsilon}\exp\!\left(\!\frac{\psi}{2\epsilon}\!\right)\!\bigg\{\!\!-\epsilon \langle\bm{D}(t),\hess\left(\psi\right)\rangle - \frac{1}{2}\lVert \bm{B}(t)^{\top}\nabla \psi \rVert_{2}^{2} \!\bigg\}.
\label{LaplacianphiANDLaplacianpsi}
 \end{align}
From (\ref{HopfColevarphi}), taking partial derivative of $\varphi$ w.r.t. $t$ gives 
{\small{\begin{equation*}
    \begin{aligned}
         &\frac{\partial \varphi }{\partial t} = \frac{1}{2\epsilon}\exp\left(\frac{\psi}{2\epsilon}\right) \frac{\partial \psi }{\partial t} \\
         						&\stackrel{\text{(\ref{HJBgenNonlinear})}}{=}  \frac{1}{2\epsilon}\exp\left(\frac{\psi}{2\epsilon}\right) \bigg\{\!\!- \frac{1}{2}\lVert \bm{B}(t)^{\top}\nabla \psi \rVert_{2}^{2} - \langle\nabla \psi, \bm{f}\rangle  - \epsilon \langle\bm{D}(t),\hess\left(\psi\right)\rangle \!\!\bigg\}  &  \\ 
         & \stackrel{\text{(\ref{nablaphiANDnablapsi}),(\ref{LaplacianphiANDLaplacianpsi})}}{=}  - \langle\nabla \varphi, \bm{f}\rangle - \epsilon\langle\bm{D}(t),\hess\left(\varphi\right)\rangle =: L_{\rm{BK}}\varphi, & 
    \end{aligned}
\end{equation*}}}
i.e., $\varphi(\bx,t)$ satisfies the backward Kolmogorov equation (\ref{BK}). 

To demonstrate that $\hat{\varphi}(\bx,t)$ satisfies the forward Kolmogorov equation (\ref{FPK}), from (\ref{HopfColevarphihat}) we compute
\begin{align}
\nabla\hat{\varphi} = \exp\left(-\frac{\psi}{2\epsilon}\right)\left(-\frac{\rho^{\text{opt}}}{2\epsilon}\nabla\psi + \nabla\rho^{\text{opt}}\right).
\label{gradphihat}
\end{align}
Consequently, 
\begin{align}
&\langle\bm{D}(t),\hess\left(\hat{\varphi}\right)\rangle = \displaystyle\sum_{i,j=1}^{n}(\bm{D}(t))_{ij}\displaystyle\frac{\partial^{2}\hat{\varphi}}{\partial x_{i}x_{j}}\nonumber\\
\stackrel{\text{(\ref{HopfColevarphihat})}}{=}& \displaystyle\sum_{i,j=1}^{n}(\bm{D}(t))_{ij}\bigg\{\!\exp\left(-\frac{\psi}{2\epsilon}\right)\frac{\partial^{2}\rho^{\text{opt}}}{\partial x_{i}\partial x_{j}} - \frac{1}{2\epsilon}\frac{\partial\rho^{\text{opt}}}{\partial x_{j}}\exp\left(-\frac{\psi}{2\epsilon}\right)\frac{\partial\psi}{\partial x_{i}} \nonumber\\
&\quad- \frac{\partial^{2}\psi}{\partial x_{i}\partial x_{j}}\frac{\rho^{\text{opt}}}{2\epsilon}\exp\left(-\frac{\psi}{2\epsilon}\right) -\frac{\partial\psi}{\partial x_{j}}\left(\frac{1}{2\epsilon}\frac{\partial\rho^{\text{opt}}}{\partial x_{i}}\exp\left(-\frac{\psi}{2\epsilon}\right) \right.\nonumber\\
&\quad\left.+ \frac{\rho^{\text{opt}}}{4\epsilon^2}\exp\left(-\frac{\psi}{2\epsilon}\right)\frac{\partial\psi}{\partial x_{i}}\right) \!\bigg\}\nonumber\\
=& \exp\left(-\frac{\psi}{2\epsilon}\right)\bigg\{\! \epsilon\langle\bm{D}(t),\hess(\rho^{\text{opt}})\rangle - \langle\nabla\rho^{\text{opt}},\bm{D}(t)\nabla\psi\rangle \nonumber\\
&\quad- \frac{\rho^{\text{opt}}}{2}\langle\bm{D}(t),\hess(\psi)\rangle + \frac{\rho^{\text{opt}}}{4\epsilon}\langle\nabla\psi,\bm{D}(t)\nabla\psi\rangle\!\bigg\}.
\label{Laplacianphihat}
\end{align}
Also, we have
\begin{align}
-\nabla\cdot\left(\rho^{\text{opt}}\bm{D}(t)\nabla\psi\right) = -\langle\nabla\rho^{\text{opt}},\bm{D}(t)\nabla\psi\rangle - \rho^{\text{opt}}\langle\bm{D}(t),\hess(\psi)\rangle,
\label{IntermedChainRule}	
\end{align}
and
\begin{align} 
         \nabla \cdot  (\hat{\varphi} \bm{f}) &= \langle\nabla \hat{\varphi},\bm{f}\rangle + \hat{\varphi} \nabla\cdot\bm{f}  \notag \\
         & \stackrel{\text{(\ref{HopfCole}),(\ref{gradphihat})}}{=}  \exp\left(-\frac{\psi}{2\epsilon}\right)\bigg(- \frac{\rho^{\text{opt}}}{2 \epsilon} \langle\nabla\psi,\bm{f}\rangle + \langle\nabla\rho^{\text{opt}},\bm{f}\rangle \nonumber\\
         &\qquad\qquad\qquad\qquad\qquad\qquad\qquad+ \rho^{\text{opt}} \nabla\cdot\bm{f}\bigg),
\label{divergenceofphihattimesf}         
\end{align} 
wherein we used the chain rule for the divergence operator. Taking partial derivative of (\ref{HopfColevarphihat}) w.r.t. $t$, we then get 
\begin{equation}
    \begin{aligned}
        \frac{\partial \hat{\varphi}}{\partial t} &= \exp\left(-\frac{\psi}{2\epsilon}\right)\bigg( \frac{\partial \rho^{\text{opt}} }{\partial t} - \frac{\rho^{\text{opt}}}{2\epsilon} \frac{\partial \psi}{\partial t}\bigg) \\
                                             &\stackrel{\text{(\ref{SBwDriftOptConds})}}{=}  \exp\left(-\frac{\psi}{2\epsilon}\right)\left[\bigg\{-\nabla\cdot(\rho^{\text{opt}}\bm{f}) - \nabla\cdot\left(\rho\bm{D}(t)\nabla\psi\right) \right.\\
                                             &\:+ \left.\epsilon\langle\bm{D}(t),\hess\left(\rho^{\text{opt}}\right)\rangle\bigg\} -\frac{\rho^{\text{opt}}}{2\epsilon}\bigg\{\! -\frac{1}{2}\parallel \bm{B}(t)^{\!\top}\nabla\psi\parallel_{2}^{2} \right.\\
                                             &\:\left.-\langle\nabla\psi,\bm{f}\rangle - \epsilon\langle\bm{D}(t),\hess(\psi)\rangle \!\bigg\}\right]\\                          
                                            & \stackrel{\text{(\ref{Laplacianphihat}),(\ref{IntermedChainRule}),(\ref{divergenceofphihattimesf})}}{=} - \nabla \cdot(\hat{\varphi}\bm{f}) + \epsilon \Delta \hat{\varphi},
                                             \end{aligned} 
\end{equation}
i.e., $\hat{\varphi}(\bx,t)$ satisfies (\ref{FPK}). 

Combining (\ref{HopfCole}) with (\ref{SBwDriftBC}), we obtain the boundary condition (\ref{FPKBKbcs}). Finally, combining (\ref{HopfColevarphi}) with (\ref{Provinguoptequalsgradpsi}), yields the optimal control $\bm{u}^{\text{opt}}(\bx,t) = 2\epsilon\bm{B}(t)^{\top}\nabla\log\varphi$.\qedsymbol

 
\vspace*{-0.1in}
\subsection{Proof of Theorem \ref{GradTVP2IVPTheorem}}\label{AppendixProofGradTVP2IVPTheorem}
\vspace*{-0.1in}
\noindent First, notice that (\ref{IVP1}) follows from combining $s=1-t$ and (\ref{BKgradrift}). Next, using 
\begin{align}
p(\bx,s) := q(\bx,s) \exp\left(-V(\bx)/\epsilon\right),
\label{pintermsofq}
\end{align}
we find that
\begin{equation}
\begin{aligned} \label{gradq}
        \nabla q &= \exp\left(V/\epsilon\right)\left( \nabla p +\frac{p}{\epsilon} \nabla V \right).  
        \end{aligned}
\end{equation}
Applying divergence operator (w.r.t. $\bx$) to both sides of (\ref{gradq}) yields
\begin{align} \label{laplaceq}
        \Delta q & = \bigg\langle \nabla p + \frac{p}{\epsilon} \nabla V, \nabla \exp\left(V/\epsilon\right)\bigg\rangle + \exp\left(V/\epsilon\right) \nabla \cdot \left( \nabla p + \frac{p}{\epsilon} \nabla V \right) \nonumber\\ 
        & = \exp\left(V/\epsilon\right) \bigg(\frac{p}{\epsilon^2} \lVert \nabla V \rVert_{2} ^{2} + \Delta p + \frac{2}{\epsilon} \langle\nabla p, \nabla V\rangle +\frac{p}{\epsilon}\Delta V \bigg).
\end{align}
From (\ref{pintermsofq}), we have
\begin{align}
\dfrac{\partial p}{\partial s} &= \exp\left(-V/\epsilon\right)\dfrac{\partial q}{\partial s}\nonumber\\
&\stackrel{\text{(\ref{IVP1})}}{=} \exp\left(-V/\epsilon\right)\left(-\langle\nabla q, \nabla V\rangle + \epsilon\Delta q\right)\nonumber\\
&\stackrel{\text{(\ref{gradq}),(\ref{laplaceq})}}{=} - \bigg\langle\nabla p +\frac{p}{\epsilon} \nabla V, \nabla V\bigg\rangle + \epsilon \left(\frac{p}{\epsilon^2} \lVert \nabla V \rVert_{2} ^{2} + \Delta p \right.\nonumber\\
&\qquad\quad\left.+ \frac{2}{\epsilon} \langle\nabla p, \nabla V\rangle +\frac{p}{\epsilon}\Delta V\right)\nonumber\\
&= \langle\nabla p, \nabla V\rangle + \epsilon \Delta p + p \Delta V = \nabla\cdot\left(p\nabla V\right) + \epsilon \Delta p,
\label{DerivingPDEinq}
\end{align}
which is indeed the PDE in (\ref{IVP2}). Setting $s=0$ in (\ref{pintermsofq}) recovers the initial condition in (\ref{IVP2}). This completes the proof. \qedsymbol

\vspace*{-0.1in}
\subsection{Proof of Theorem \ref{DegenTVP2IVPTheorem}} \label{AppendinxDegenTVP2IVPTheorem}
\vspace*{-0.1in} 
\noindent Notice that (\ref{BKdegenIVP}) follows by combining $s:=1-t$, and (\ref{BKdegen}). To derive (\ref{FKdegenIVP}), we start by taking the gradient of
\begin{align} \label{qprelation}
\widetilde{p}(\bm{\xi},-\bm{\eta},s)=q(\bm{\xi},\bm{\eta},s)\!\exp\left(\! -\frac{1}{\epsilon} \left(\!\frac{1}{2} \lVert \bm{\eta}\rVert_{2}^{2} + V(\bm{\xi}) \!\right)\!\right)
\end{align}
w.r.t. $\bm{\xi}$ and $\bm{\eta}$, respectively, to obtain
\begin{subequations}
\begin{align}
\nabla_{\bm{\xi}} q &=\exp\left( \frac{1}{\epsilon} \left( \frac{1}{2} \lVert \bm{\eta}\rVert_{2}^{2} + V(\bm{\xi}) \right) \right)
 \left(\nabla_{\bm{\xi}} \widetilde{p} + \frac{\widetilde{p}}{\epsilon} \nabla_{\bm{\xi}}  V  \right), \label{xigrad} \\
  \nabla_{\bm{\eta}} q &=\exp\left( \frac{1}{\epsilon} \left( \frac{1}{2} \lVert \bm{\eta}\rVert_{2}^{2} + V(\bm{\xi}) \right) \right) \left(-\nabla_{-\bm{\eta}} \widetilde{p} + \frac{1}{\epsilon} \widetilde{p} \:\bm{\eta} \right). \label{etagrad}
 \end{align}
 \label{Gradients}
\end{subequations}
Applying divergence operator w.r.t. $\bm{\eta}$ on both sides of (\ref{etagrad}) yields 
\begin{align} \label{etalaplace}
\Delta_{\bm{\eta}} q &= \exp\left( \frac{1}{\epsilon} \left( \frac{1}{2} \lVert \bm{\eta}\rVert_{2}^{2} + V(\bm{\xi}) \right) \right) \left(  \Delta_{-\bm{\eta}} \widetilde{p} -\frac{2}{\epsilon} \langle \nabla_{-\bm{\eta}}\widetilde{p},\bm{\eta} \rangle \right. \nonumber \\
&+\left. \frac{\widetilde{p}}{\epsilon}+ \frac{\widetilde{p}}{\epsilon} \lVert \bm{\eta} \rVert_{2}^{2}\right).
\end{align}
Thus, we have 
\begin{align*}
&\frac{\partial p}{ \partial s} =  \frac{\partial \widetilde{p}}{\partial s} \stackrel{\text{(\ref{qprelation})}}{=} \exp\left( -\frac{1}{\epsilon} \left( \frac{1}{2} \lVert \bm{\eta}\rVert_{2}^{2} + V(\bm{\xi}) \right) \right) \frac{\partial q} {\partial s}  \\
&\stackrel{\text{(\ref{BKdegenIVP}),(\ref{Gradients}),(\ref{etalaplace}) }}{=} \exp\left( -\frac{1}{\epsilon} \left( \frac{1}{2} \lVert \bm{\eta}\rVert_{2}^{2} + V(\bm{\xi}) \right) \right)  \left\{ \langle \bm{\eta}, \nabla_{\bm{\xi}} \widetilde{p} + \frac{\widetilde{p}}{\epsilon} \nabla_{\bm{\xi}}  V   \rangle  \right. \nonumber \\
&\qquad\qquad -\langle \nabla_{\bm{\xi}}V\left(\bm{\xi}\right) - \kappa\bm{\eta}, -\nabla_{-\bm{\eta}} \widetilde{p} + \frac{1}{\epsilon} \widetilde{p}  \bm{\eta}\rangle \nonumber  \\
& \qquad\qquad\left. + \Delta_{-\bm{\eta}} \widetilde{p} -\frac{2}{\epsilon} \langle \nabla_{-\bm{\eta}}\widetilde{p},\bm{\eta} \rangle 
 \frac{\widetilde{p}}{\epsilon}+ \frac{\widetilde{p}}{\epsilon} \lVert \bm{\eta} \rVert_{2}^{2}\right\} \nonumber \\ 
 &= \exp\left( -\frac{1}{\epsilon} \left( \frac{1}{2} \lVert \bm{\eta}\rVert_{2}^{2} + V(\bm{\xi}) \right) \right) \bigg\{  \langle \bm{\eta},\nabla_{\bm{\xi}} \widetilde{p} \rangle -\kappa \langle \nabla_{-\bm{\eta}}\widetilde{p},\eta \rangle  \nonumber \\
 &\quad + \langle \nabla_{-\bm{\eta}}\widetilde{p},\nabla_{\bm{\xi}} V\rangle  + \kappa \widetilde{p} + \epsilon \kappa \Delta_{-\bm{\eta}} \widetilde{p}  \bigg \} \nonumber \\
 &= 
 \exp\left( -\frac{1}{\epsilon} \left( \frac{1}{2} \lVert \bm{\vartheta}\rVert_{2}^{2} + V(\bm{\xi}) \right) \right) \bigg\{  \langle -\bm{\vartheta},\nabla_{\bm{\xi}} p \rangle + \kappa \langle \nabla_{\bm{\vartheta}}p,\vartheta \rangle  \nonumber \\
 &\quad + \langle \nabla_{\bm{\vartheta}}p,\nabla_{\bm{\xi}} V\rangle  + \kappa p + \epsilon \kappa \Delta_{\bm{\vartheta}} p  \bigg \} \nonumber \\
  & =
 \exp\left( -\frac{1}{\epsilon} \left( \frac{1}{2} \lVert \bm{\vartheta}\rVert_{2}^{2} + V(\bm{\xi}) \right) \right) \bigg\{
 -\langle \bm{\vartheta},\nabla _{\bm{\xi}}p  \rangle \nonumber \\
 &\quad + \nabla_{\bm{\vartheta}} \cdot\left(p \left(\nabla_{\bm{\xi}}V\left(\bm{\xi}\right) + \kappa\bm{\vartheta}\right)\right) + \epsilon \kappa \Delta_{\bm{\vartheta}} p \bigg \},
\end{align*}
which is the PDE  in (\ref{FKdegenIVP}). Setting $s=0$ in (\ref{qprelation}) recovers the initial condition in (\ref{FKdegenIVP}) .\qedsymbol

\vspace*{-0.1in}
\subsection{Regularity of the Transition Densities for (\ref{transPDE}) and (\ref{transPDEdegen})} \label{AppendinxRegularity} 
\vspace*{-0.1in}
\noindent In this Appendix, we point out that the transition probability densities $K(s,\bm{y},t,\bx)$ for (\ref{transPDE}) and (\ref{transPDEdegen}), are indeed positive and continuous in $\bx,\by\in\mathbb{R}^{n}$ for all $s<t$. First, recall that the transition densities themselves solve the same PDEs as in (\ref{transPDE}) and (\ref{transPDEdegen}), with initial condition $\lim_{t\downarrow s}K(s,\bm{y},t,\bx)=\delta(\bm{x}-\bm{y})$. From the maximum principle for parabolic PDEs, it follows that the transient solutions of (\ref{transPDE}) and (\ref{transPDEdegen}) are positive as long as the initial conditions are positive. The continuity, in the gradient drift case (\ref{transPDE}), is standard assuming $V\in C^{2}\left(\mathbb{R}^{n}\right)$; see e.g., \cite[Ch. 1.2]{stroock2008partial}. For the degenerate diffusion case (\ref{transPDEdegen}), the situation is more subtle: a result from Villani \cite[Theorem 7]{villani2006hypocoercivity} ensures the continuity (w.r.t. both state and time) of the transient solutions, under the assumptions on $V(\cdot)$ stated in Section \ref{SubsectionSchrodingerSystemwithMixedDrift}, viz. $V\in C^{2}\left(\mathbb{R}^{m}\right)$, $\inf V > -\infty$, and $\hess\left(V\right)$ uniformly lower bounded. 


\bibliographystyle{IEEEtran}
\bibliography{references}

\end{document}